\documentclass[11pt]{article}%
\usepackage[latin1]{inputenc}
\usepackage[T1]{fontenc}
\usepackage{amsfonts,amssymb,amsmath, epsfig}
\usepackage{color,graphicx,graphics}
\usepackage{amsmath,amstext,amssymb,amsfonts, amscd}
\usepackage{hyperref}
\usepackage{a4wide}
\usepackage{booktabs}
\usepackage{cite}

\usepackage{cleveref}
\numberwithin{equation}{section}
\def\Box{\leavevmode\vbox{\hrule
		\hbox{\vrule\kern4pt\vbox{\kern4pt}%
			\vrule}\hrule}}

\newcounter{appendix}
\def\appendix{\advance\c@appendix by 1
	\def\thesection{\Alph{section}}
	\ifnum\c@appendix=1 \setcounter{section}{-1} \fi
	\@startsection {section}{1}{\z@}{-3.5ex plus -1ex minus
		-.2ex}{2.3ex plus .2ex}{\Large\bf}}

\newtheorem{lemma}{Lemma}[section]

\newtheorem{remark}{Remark}[section]
\def\softd{{\leavevmode\setbox1=\hbox{d}%
          \hbox to 1.05\wd1{d\kern-0.4ex{\char039}\hss}}}

\title{Solving Random Hyperbolic Conservation Laws Using Linear Programming }
\author{Shaoshuai Chu,\footnote{RWTH Aachen University, Institut f\"ur Geometrie und Praktische Mathematik, 52056 Aachen, Germany, \\{\tt \{chu, herty, zhou\}@igpm.rwth-aachen.de}}\quad  Michael Herty,$^*$~ 
 M\'aria Luk\'a\v{c}ov\'a-Medvi{\softd}ov\'a,\footnote{Johannes Gutenberg University Mainz, Institute of Mathematics, Staudingerweg 9, 55128 Mainz, Germany, {\tt lukacova@uni-mainz.de}}~~ and~ Yizhou Zhou$^*$}
\begin {document}
\date{}
\maketitle
\begin{abstract}
A novel structure-preserving numerical method to solve random hyperbolic systems of conservation laws is presented. The method uses a concept of generalized, measure-valued solutions to random conservation laws.  This yields a linear partial differential equation with respect to the Young measure and allows to compute the approximation based on linear programming problems. We analyze structure-preserving properties of the derived numerical method and discuss its advantages and disadvantages. We numerically demonstrate the approach on the one-dimensional Burgers and isentropic Euler equations and compare with stochastic collocation. In addition, we introduce a discontinuous-flux test in which different entropies used in the linear-program objective select different weak entropy solutions, and we report the corresponding changes in the moments and supports of the Young measure.
\end{abstract}



\section {Introduction}
Random nonlinear conservation laws arise in various engineering, atmospheric or geophysical applications. In recent years, several numerical methods have been proposed for random conservation laws; see, e.g., \cite{MR3674794} for an overview. Popular methods include non-intrusive (Quasi-)Monte Carlo simulations \cite{MR3674794,MR3793404,MR3471175}, stochastic collocation \cite{MR2199923,MR2723020,MR3202521,MR4045233}, or (intrusive) generalized polynomial chaos expansion \cite{MR3328389,MR2723020,MR2501693,MR4066974,MR1507356} and references therein.  

The advantages and disadvantages of different approaches have been discussed extensively in the literature and we refer to  \cite{chertockchallenges} for recent discussions. At this point, we only mention that the Monte Carlo methods, while robust, are not necessarily the most efficient. Similar to the stochastic collocation method, the numerical methods provided for the deterministic case can be easily reused, strongly simplifying their applicability. In contrast, intrusive methods, based on generalized polynomial chaos expansion, known as the stochastic Galerkin methods, are less flexible. They require typically projection of the Galerkin approximation  that leads to a  new deterministic system of equations for the coefficients of the expansion~\cite{MR2660316,MR3129545}. Provided the solution is sufficiently smooth, this intrusive approach is usually more accurate even for a low degree of the expansion \cite{MR2823001,MR2855645}.  However, the resulting system for the coefficients may be challenging to solve due to a possible loss of hyperbolicity \cite{MR3202523,MR4228317,MR3968112,MR4066974,MR4772617,MR4699586}. Moreover,  it also requires a re-implementation of the numerical method.  The  stochastic collocation as well as stochastic Galerkin methods may encounter Gibbs phenomena  across jump discontinuities due to the spectral approximation requiring additional filtering \cite{MR1874071,MR4045233,MR4228317,MR2333926,MR3874527}, spectral viscosity \cite{MR977947,MR1770237}, or adaptive diffusion \cite{chertockchallenges,MR2979844} techniques.

In this work, a different way to numerically compute approximate solutions to random nonlinear hyperbolic partial differential equations (PDEs)  \eqref{PDE1} is presented. The main idea is to consider a solution in a measure-valued framework as parameterized Young measures instead of the weak (distributional) solutions \cite{MR3409135,MR257325} (see also \cite{DiPerna1985}, where the concept of measure-valued solutions for conservation laws was introduced in a rigorous analytical framework and has inspired numerous developments). This yields a linear problem that allows for a different treatment circumventing  most of the mentioned problems. Our approach is motivated  by the recent works on dissipative measure-valued solutions to the compressible fluid flows, see e.g., \cite{MR4390192,MR4130543,MR4674068}. We propose here a novel numerical method based on the measure-valued formulation of a hyperbolic conservation law. In the aforementioned publications the Young measure is generated by a weakly convergent approximated sequence and represents the underlying (dissipative) measure-valued solution  of the compressible Euler and Navier-Stokes equations.  Contrary to this approach, in the present work  the Young measure is directly computed in the construction of the numerical scheme. Our new approach is presented here for general one-dimensional systems of conservation laws and numerical results are given for the isentropic Euler system, the Burgers equation, and a discontinuous flux test, while in \cite{MR4390192,MR4130543,MR4674068} the multidimensional fluid dynamic equations are considered. We also refer to a related recent publication \cite{MR4777936}, which proposes a moment-based framework for parameter-dependent hyperbolic conservation laws. Their approach exploits the linearity of the moment equations and relies on the Moment-SOS hierarchy to approximate entropy solutions, which allows for the proof of convergence under certain structural assumptions.

The difference is that there the moments are computed using a sufficiently large family of entropies of the system. The latter is, however, not necessarily available for a general hyperbolic system. Furthermore, also moments with respect to time and space are considered, whereas we propose to use only moments in the random dimension. This yields a Young measure parameterized by $(t,x)$ and allows to use standard finite-volume schemes for the propagation of the moments of the measure.  Note that there is a similarity with \cite{MR3828934,MR4045233}, where a moment-based approach to random conservation laws has been proposed. However, in the present paper we provide a different formulation of the closure problem compared to the above references. Furthermore, the presented approach does not require filters \cite{MR4045233} and it is applicable also to the system case. This is at the expense of solving a linear programming problem in contrast  to the nonlinear optimization problem in the aforementioned publications. As far as we are aware the present methodology is novel and  there are no further similarities with existing works. The remaining part of the paper is organized as follows. In the next section, the method will be formally derived and some structure-preserving properties will be discussed.  Results of numerical simulations for the random Burgers equation, the isentropic Euler system, and one example taking the discontinuous flux will be given in \S\ref{sec:num}.  A discussion and outlook  section closes this work.

%%%%%%%%%%%%%%%%%%%%%%%%%%%%%%%%%%%%%%%%%%%%%%%%%%%%%%%%%%
\section{The Method}
Let $\big(\Omega, \mathcal{F}(\Omega), \mathbb{P} \big)$ be a  probability space with $\Omega$ being  a metric space of events, $\mathcal{F}(\Omega)$ $\sigma$-field of Borel subsets of $\Omega$, and $\mathbb{P}$ a probability measure on $\Omega.$ 

We consider the random Cauchy problem for a system of nonlinear hyperbolic conservation laws 
\begin{eqnarray}
\label{PDE1}
&& \partial_t u(t,x,\xi) + \partial_x f(u(t,x,\xi))  = 0,   \qquad   t \in (0,T),\,\,\, x \in \mathbb{R}\,\,\, \mbox{ and }\, \xi \in \Omega,\nonumber \\ 
&& u(0,x,\xi)=u_0(x,\xi) \in \Bbb{R}^n, \qquad  \qquad  \  \   x \in \mathbb{R}\,\, \mbox{ and }\,\, \xi \in \Omega.
\end{eqnarray}
For simplicity, the randomness enters only through the random initial data $u_0$ that are Borel measurable mapping from $\Omega$ to the phase space in $\mathbb{ R}^n.$

The nonlinear flux $f : \mathbb{R}^n \rightarrow \mathbb{R}^n$ is assumed to be hyperbolic, i.e., the Jacobian matrix $A=\frac{\partial f}{\partial u}$ has only real eigenvalues and is diagonalizable. For a detailed definition and discussion of hyperbolicity, we refer to \cite[\S 3]{MR3468916}. We also assume that \eqref{PDE1} is accompanied by a convex entropy $\eta$ and the corresponding entropy flux pair $q$ such that the entropy inequality 
\begin{align}
\label{PDE2}
\partial_t \eta(u(t,x,\xi)) + \partial_x q(u(t,x,\xi))  \leq  0
\end{align}
holds true. In what follows we consider a very weak, i.e., measure-valued, formulation of (\ref{PDE1})-(\ref{PDE2}), cf.~\eqref{PDE-YM}, \eqref{PDE-YM-entropy}.   Existence of measure-valued solutions, their compatibility and measure-valued-strong uniqueness principle for multidimensional hyperbolic conservation laws \eqref{PDE1}-\eqref{PDE2}  were investigated in \cite{MR4390192}. 

For theoretical results on existence and uniqueness of weak entropy solutions for (deterministic) one-dimensional systems of hyperbolic conservation laws and multidimensional scalar hyperbolic conservation laws we refer to \cite{MR4777936,MR1816648,MR3468916,MR3793404}.

The situation changes dramatically for multidimensional systems of hyperbolic conservation laws. As shown in the pioneering works of 
De Lellis and Sz\'ekelyhidi~\cite{DLSz, DlSz1}   many weak entropy solutions may exist,  see also \cite{Klingenberg1} and the references therein. As shown recently in \cite{CF} the set of initial data for which infinitely many weak entropy solutions of the compressible Euler equations exists is dense in $L^p$-topology. We also refer to recent works, where various selection criteria were studied to recover well-posedness of the compressible Euler equations in the class of measure-valued or dissipative (measure-valued) solutions \cite{Breit, klingenberg2025, FJL, feireisl-lukacova2025}.

The challenges of well-posedness in the context of parameter-dependent (random) hyperbolic systems are addressed in \cite{MR4777936}, which shows that existence and uniqueness can be nontrivial and require careful moment-based formulations and closure assumptions. This underlines the complexity of extending deterministic theories to the random hyperbolic problems.  In particular, the random conservation law \eqref{PDE1} can also be understood as a deterministic conservation law on the extended phase space $(x,\xi)$ using properly weighted solution spaces, see, e.g., \cite{MR4699586}.  \\
 
We proceed with the formal derivation of the numerical approximation of the Cauchy problem \eqref{PDE1}-\eqref{PDE2}. Without loss of generality we assume that the probability distribution is an absolutely continuous random variable on $\Omega \subseteq \mathbb{R}$ with respect to the Lebesgue measure. Hence, there exists an essentially bounded probability density $p:\Omega \to [0,\infty).$ Further, we consider the set  $\{ \phi_i \}_{i=1}^\infty$ of  orthonormal polynomials on $\Omega$ with values in $\mathbb{R}$ and  of degree at most $i,$ such that 
\begin{equation*}
\int_{\Omega} \phi_i(\xi) \; \phi_j(\xi)  \; p(\xi) \mathrm{d}\xi = \delta_{i,j}, \; i,j=1, \dots
\end{equation*}
Here, $\delta_{i,j}$ denotes the Kronecker symbol. For uniform and Gaussian distributions, different orthogonal polynomials are suitable; e.g., Legendre and Hermite, respectively. A comprehensive list can be found in \cite{MR2723020,MR372517}.  Flexibility in the choice of suitable basis functions  allows to present our new method in a uniform framework. 
Formally,   the case of stochastic collocation approach in the aforementioned formulation is obtained by setting $\phi_i(\xi) = \delta(\xi-\xi_i)$ for a sequence of points $\{ \xi_i \}_i$ with $i=1, 2,3,\dots.$ In \S\ref{sec:num} we present results of numerical simulations obtained by using these basis functions. 
\begin{remark}
 The stochastic collocation approach is a non-intrusive method to solve random PDEs. It approximates the solution by evaluating deterministic PDE solvers at selected collocation points in the random space. Formally, this corresponds to choosing basis functions of the form $\varphi_i(\xi) = \delta(\xi - \xi_i)$, where $\{\xi_i\}$ are collocation points. For each $\xi_i$, a deterministic version of the PDE is solved independently, and statistical quantities (such as mean and variance) are computed via numerical quadrature. This method is particularly attractive because it allows the reuse of existing deterministic solvers without modification. Moreover, it handles discontinuities well and avoids the Gibbs phenomenon associated with spectral methods. However, it may require a large number of collocation points in high-dimensional random spaces to achieve sufficient accuracy. For more details, we refer the readers to \cite{MR2199923,MR2723020,MR3202521,MR4045233}.
 \end{remark}
 
 \begin{remark}
 It is important to note that in \cite{MR4777936} the main point is to have many entropies available, since those allow for the reconstruction of the measure. This is certainly the case in the scalar case, but may not be true for systems. Here, we do not rely on those moments for the reconstruction. The moments we have used are only with respect to the variable $\xi$ (not $u$).
 \end{remark}

\subsection{Minimization Formulation for Stochastic Galerkin Approach}
Recall that in the stochastic Galerkin approach, we multiply equation \eqref{PDE1} by $p(\xi)$,  and by $\phi_i(\xi)$. Integrating over $\Omega$ yields an equation for the $i$-th stochastic moment $U_i$ of $u$ given by  
\begin{align} 
\label{Galerkin moments} U_i(t,x) := \int_{\Omega} \phi_i(\xi) u(t,x,\xi) p(\xi) \mathrm{d}\xi, 
\end{align}
where $U_i$ fulfills the PDE 
\begin{align}\label{pde moment}
\partial_t U_i(t,x) + \partial_x \int_{\Omega} \phi_i(\xi) f(u(t,x,\xi)) p(\xi) \mathrm{d}\xi = 0, \qquad i=1,2,3,\dots
\end{align}
The initial data for {the evolution of $U_i$} is obtained by the projection of $u_0,$ i.e., 
$$
U_{i}^0(x) = \int_{\Omega} \phi_i(\xi) u_0(x,\xi) p(\xi) \mathrm{d}\xi, \qquad i=1,2,3,\dots
$$
Note that in the case of stochastic collocation {approach} it holds $U_i(t,x) = u(t,x,\xi_i)$ for  given collocation points $\xi_i \in \Omega.$ The reconstruction of the flux term $f(u(t,x,\xi))$ out of the moments $U_i$ is the major task. In the stochastic Galerkin framework, e.g., a (truncated) series expansion $u(t,x,\xi) = \sum_{i=1}^{N_\xi} U_i(t,x) \phi_i(\xi)$ is used, requiring a projection of $f(u)$ on the space generated by $\{\phi_i \}_{i=1}^{N_\xi}.$  In the case of truncated series expansions this procedure  may lead to loss of hyperbolicity depending on the formulation of expansion variables, see \cite{MR2501693,MR4066974}. 

Inspired by kinetic theory a possibility to determine $u=u(t,x,\xi)$ is to consider a minimization problem with constraints given by the moments $U_i$;  see \cite{MR3828934,MR2501693}. 

For each fixed $(t,x)$ and given moments $U_i(t,x)$ for $i=1,\dots,N_\xi$, a function $\xi \to u^*(\xi; t,x )$ is obtained as a solution to problem \eqref{moment nonlinear}.
\begin{align}\label{moment nonlinear}
\hspace{-0.3cm}u^* =\mbox{argmin}_{v(\cdot)}  \int_{\Omega} \eta( v(\xi) ) p(\xi) \mathrm{d}\xi\,\,\, \mbox{ subject to} \int_{\Omega} \phi_i(\xi) v(\xi) p(\xi) \mathrm{d}\xi = u_i(t,x),  
\end{align}
for $i=1,\dots,N_\xi.$
Then equation \eqref{pde moment} is closed in terms of $U_i$ by approximating  
\begin{equation*} 
\int_{\Omega} \phi_i(\xi) f(u(t,x,\xi)) p(\xi) \mathrm{d}\xi \approx \int_{\Omega} \phi_i(\xi) f(u^*(\xi; t,x)) p(\xi) \mathrm{d}\xi.
\end{equation*}
Namely, we take $\int_{\Omega} \int_{\Omega} \phi_i(\xi) f(u^*(\xi; t,x)) p(\xi) \mathrm{d}\xi$ as the flux term in computing the equation \eqref{pde moment}. Here $u^*$ is the solution to the problem \eqref{moment nonlinear}.
More details on this approach are given, e.g., in \cite{MR3828934}. Since problem \eqref{pde moment} is nonlinear, its solution is non-trivial and it is possibly computationally expensive. Our scheme proposed here is motivated by a similar minimization problem but uses a formulation in terms of the Young measures to introduce an additional linear structure. 

\subsection{Formal Derivation of the Proposed Scheme}
\renewcommand{\u}{\mathbf{u}}
The proposed method  aims to approximate the stochastic moments $U_i$  of the random solution $u(t,x,\xi)$ as in equation \eqref{Galerkin moments}. We use a linear structure of the measure-valued formulation of \eqref{PDE1}-\eqref{PDE2}.  Denote by $\mathcal{P}(\mathbb{R}^n)$ the space of probability measures on $\mathbb{R}^n.$

We say that the Young measure $\nu_{t,x,\xi} \in \mathcal{P}(\mathbb{R}^n)$ is a measure-valued solution of \eqref{PDE1}-\eqref{PDE2} if 
\begin{align}\label{PDE-YM}
\int_{0}^T \int_{\Bbb R} \int_\Omega \partial_t \psi(t,x) \phi(\xi) \int_{\mathbb{R}^n} \u \;  \mathrm{d}\nu_{t,x,\xi}(\u)   + \partial_x  \psi(t,x) \phi(\xi) \int_{\mathbb{R}^n} f(\u) \; \mathrm{d}\nu_{t,x,\xi}(\u) p(\xi)  \mathrm{d} \xi \mathrm{d} x \mathrm{d} t  = 0,
\end{align}
for all test functions $\psi \in C^1_c((0,T) \times \Bbb R),$ $\phi \in C(\Omega) $ and
\begin{align}\label{initial-ym-0}
\nu_{0,x,\xi}(\cdot) = \delta(\cdot - u_0(x,\xi) ) \mbox{ for  a.a. } x  \in \Bbb R, \xi \in \Omega.
\end{align}
Further, we require that the entropy inequality holds in the following sense
\begin{align} \label{PDE-YM-entropy}
\int_{0}^T \int_{\Bbb R} \int_\Omega \partial_t \tilde{\psi}(t,x) \tilde{\phi}(\xi) \int_{\mathbb{R}^n} \eta(\u) \;  \mathrm{d}\nu_{t,x,\xi}(\u)   + \partial_x  \tilde{\psi}(t,x) \tilde{\phi}(\xi) \int_{\mathbb{R}^n} q(\u) \; \mathrm{d}\nu_{t,x,\xi}(\u) p(\xi)  \mathrm{d} \xi \mathrm{d} x \mathrm{d} t  \geq 0,
\end{align}
for all test functions $\tilde{\psi} \in C^1_c((0,T) \times \Bbb R),$ $\tilde{\phi} \in C(\Omega) $, $\tilde{\psi} \geq 0, \tilde{\phi} \geq 0.$

Clearly, provided that $\nu_{t,x,\xi}$ is a Dirac measure concentrated at $u=u(t,x,\xi)$ for every $t\in[0,T],x\in \mathbb{R},\xi\in\Omega$,  the standard weak formulation of \eqref{PDE1}, \eqref{PDE2} is recovered. Due to the linearity, there may be many solutions to equation \eqref{PDE-YM} since equation \eqref{PDE-YM} is linear in $\nu_{t,x,\xi}.$

Similar to the stochastic Galerkin approach presented in the previous section, we take for the test function $\phi_i$ in equation \eqref{PDE-YM}  to obtain $i=1,\dots,N_\xi$ equations for the stochastic moments $u_i(t,x)$ of the measure $\nu_{t,x,\xi},$ i.e., 
\begin{equation*}
u_i(t,x) = \int_{\Omega} \int_{\mathbb{R}^n} \phi_i(\xi)  \; \u \; \mathrm{d} \nu_{t,x,\xi}(\u) \; p(\xi)  d\xi. 
\end{equation*} 

The initial data for the stochastic moment of $u_i$ are obtained by projection using the test function $\phi_i$:
\begin{align}  \label{PDE-YM-ID} 
u_{i}^0(x) =  \int_{\Omega} \int_{\mathbb{R}^n}  \phi_i(\xi)  \; \u \; \mathrm{d} \nu_{0,x,\xi}(\u) \; p(\xi)  \mathrm{d}\xi = U_{i}^0(x), \; i=1,\dots,N_\xi.
\end{align}

The evolution equation for $u_i$ is obtained by testing \eqref{PDE-YM} with the test function $\phi_i,$ \ $i=1,\dots,N_\xi$  yielding the weak formulation of the PDE
\begin{align}\label{PDE-YM-moments}
\partial_t u_i(t,x) + \partial_x    \int_{\Omega} \int_{\mathbb{R}^n} \phi_i(\xi) f(\u) \; \mathrm{d}\nu_{t,x,\xi}(\u) p(\xi) \mathrm{d}\xi = 0.
\end{align}

Thus, the moments $u_i$ are approximations to the stochastic moments $U_i$ (see equation \eqref{Galerkin moments}) of the solution $u(t,x,\xi)$ to the original random hyperbolic conservation law \eqref{PDE1}. Similarly to \eqref{pde moment} the system \eqref{PDE-YM-moments} is not in closed form with respect to the stochastic moments $u_i$ of $\nu_{t,x,\xi.}$ Therefore,  an optimization problem, similarly to problem \eqref{moment nonlinear}, is formulated to determine a suitable measure $\nu$. The constraints are given by the stochastic  moments of $\nu_{t,x,\xi}.$ Mimicking problem \eqref{moment nonlinear} the proposed method minimizes the entropy over the space of probability measures. For { a.a.} $(t,x) \in \Bbb R^+ \times \Bbb R$  a family of parameterized measure  $\{ \mu^*_{\xi;t,x} \}_{\xi \in \Omega} $ on $\mathbb{R}^n$ is obtained such that  
\begin{equation}\label{ym lin prog}
\begin{aligned}
\mu^*_{\xi; t,x} =& \mbox{argmin}_{ \xi \in \Omega, \mu_\xi \in \mathcal{M}}  \int_{\Omega} \int_{\mathbb{R}^n}  \eta(\u) \mathrm{d} \mu_\xi(\u) p(\xi) \mathrm{d}\xi \;  \\ 
\mbox{\rm{ subject to } } &   \int_{\Omega} \int_{\mathbb{R}^n} \phi_i(\xi) \u \mathrm{d}\mu_\xi(\u) p(\xi) \mathrm{d}\xi = u_i(t,x), \; i=1,\dots,N_\xi. 
\end{aligned}
\end{equation}

Here, $\mathcal{M}\subset\mathcal{P}(\mathbb{R}^n)$ is a subset of the probability space. It is a freedom in this general framework and we will specify it later.

\begin{remark}
Here we can also use the notation as that in \cite{DiPerna1985}. Namely, we formulate the problem \eqref{ym lin prog} as 
\begin{equation*}
\begin{aligned}
\mu^*_{\xi; t,x} =& \mbox{argmin}_{ \xi \in \Omega, \mu_\xi \in \mathcal{M}}  \int_{\Omega} \langle \eta, \mu_{\xi} \rangle  p(\xi) \mathrm{d}\xi \;  \\ 
\mbox{\rm{ subject to } } &   \int_{\Omega} \langle \mathrm{id}, \mu_{\xi} \rangle \phi_i(\xi)  p(\xi) \mathrm{d}\xi = u_i(t,x), \quad i=1,\dots,N_\xi. 
\end{aligned}
\end{equation*}
with $\mathrm{id}(u) := u$ and
$$
\langle \eta, \mu_{\xi} \rangle := \int_{\mathbb{R}^n} \eta(u)\, d\mu_{\xi}(u).
$$
\end{remark}

\begin{remark}
    If the feasible set $\mathcal{M}$ is taken as the full probability space $\mathcal{P}(\mathbb{R}^n)$, then the following argument holds: if $\mu^{*}_{\cdot;t,x}$ is any minimizer, then by Jensen's inequality, the Dirac measure
$$
\tilde{\mu}^{*}_{\xi;t,x} := \delta_{u^{*}(\xi;t,x)} \quad \text{where} \quad u^{*}(\xi;t,x) := \langle \mathrm{id}, \mu^{*}_{\xi;t,x} \rangle
$$
will also be a minimizer. In fact, if $\eta$ is strictly convex, then $\tilde{\mu}^{*}$ will have strictly less entropy than $\mu^{*}$. However, if we take the feasible set $\mathcal{M}$ to be the probability space with a given support, then the Dirac measure may be excluded, e.g., see the numerical example in \S3.3.
\end{remark}

\begin{remark}
Since the moments of the measure are fixed in our method, we do not expect to capture two possible solutions, e.g., $\delta_{u_{1}(\xi;t,x)}$ and $\delta_{u_{2}(\xi;t,x)}$, simultaneously. As aforementioned, this non-uniqueness may happen for multidimensional hyperbolic systems.  
Another approach is to replace the constraints in \eqref{ym lin prog} by the equation \eqref{PDE-YM-moments}. In other words, we look for a measure $\mu$ such that it minimizes the total entropy, and satisfies the
$N_\xi$ equations of conservation \eqref{PDE-YM-moments}.
\end{remark}

Formulation \eqref{ym lin prog} yields a linear programming problem compared  to nonlinear problem \eqref{moment nonlinear}. This is due to the formulation by means of the Young measures.  We will exploit this fact further in the following paragraph. Family of solutions  $\mu^*_{\xi; t,x}$ is used to approximate the measure-valued solution $\nu_{t,x,\xi}$ by the following closure relation for each $i=1,\dots,N_\xi$
\begin{align}\label{YM closure}
\int_{\Omega} \int_{\mathbb{R}^n} \phi_i(\xi) f(\u) \; \mathrm{d}\nu_{t,x,\xi}(\u) p(\xi) \mathrm{d}\xi  \approx \int_{\Omega} \int_{\mathbb{R}^n} 
\phi_i(\xi) f(\u) \mathrm{d}\mu^*_{\xi; t,x }(\u) p(\xi) \mathrm{d}\xi. 
\end{align}
Namely, we take $\int_{\Omega} \int_{\mathbb{R}^n} 
\phi_i(\xi) f(\u) \mathrm{d}\mu^*_{\xi; t,x }(\u) p(\xi) \mathrm{d}\xi$ as the flux term in computing the equation \eqref{PDE-YM-moments}.
The evolution of the stochastic moments $u_i$ is fully determined by \eqref{PDE-YM-ID}, \eqref{PDE-YM-moments}, \eqref{ym lin prog}, and \eqref{YM closure}, respectively. The linearity of problem \eqref{ym lin prog} will translate into linear programming problems after discretization.  The solution $\mu^*_\xi$ to problem \eqref{ym lin prog} might not be unique. We will discuss below the consequences, but note that we are only interested in the resulting value of the flux after applying the closure \eqref{YM closure}.  
Each solution $\mu^*_{\xi; t,x}$ depends solely on the moments $u_i(t,x)$ with $i=1,\dots,N_\xi.$  Hence, for a.a. $(t,x)$ and $N_\xi$ given moments $\vec{u}=\left( u_1, \dots, u_{N_\xi} \right)^T$,   each minimizer to \eqref{ym lin prog} may be denoted by
\begin{align}\label{solution lin prog}
 \mu^*_\xi(\cdot; \vec{u}(t,x) ) := \mu^*_{\xi; t,x}(\cdot).
\end{align}
Due to equation \eqref{solution lin prog}, the closure relation \eqref{YM closure}, the flux in \eqref{PDE-YM-moments} is given in terms of $\{ u_i \}_{i=1}^{N_\xi}$. Hence, the weak formulation of \eqref{PDE-YM-moments} (in the sense of \eqref{PDE-YM}) is a conservation law to be discretized by a suitable numerical method, e.g., structure-preserving  finite-volume scheme. In particular, discrete formulation of the entropy inequality \eqref{PDE-YM-entropy} needs to hold.

In the subsequent section a first-order entropy stable finite-volume scheme is presented. 

\subsection{Semi-Discrete Scheme with the Lax-Friedrichs Flux}
In what follows, we work with an equidistant spatial grid $\{C_j\}_{j=1}^{N_x}$;  $C_j$  is  a mesh cell of size $\Delta x$ with center $x_j$. The cell averages are $u_{i,j}(t) \approx \frac{1}{\Delta x}  \int_{C_j} u_i(t,x) dx.$ For the initial data the cell averages are obtained as 
\begin{equation}\label{2.14aaa}
u_{i,j}^0 = \frac1{\Delta x} \int_{C_j} \int_{\Omega} \int_{\mathbb{R}^n}  \phi_i(\xi)  \; \u \; \mathrm{d} \nu_{0,x,\xi} \; p(\xi)  \mathrm{d}\xi,
\end{equation} 
see equation \eqref{PDE-YM-ID}. We approximate the weak formulation of equation \eqref{PDE-YM-moments} by means of the first-order Lax-Friedrichs finite-volume method that  yields an approximation $u^n_{i,j}$ to $u_{i,j}(t_n)$ at time $t_n=n \Delta t$ for $n=0,1,\dots$: 
\begin{align}\label{YM-LxF}
u_{i,j}^{n+1} = \frac12 \left( u^n_{i,j+1} + u^n_{i,j-1} \right) - \frac{ \Delta t}{2 \Delta x} \left( \mathcal{F}_{i,j+1}^n  - 
\mathcal{F}_{i,j-1}^n \right),\,\, i=1,\dots,N_\xi, \quad j=1,\dots,N_x.
\end{align}  
The first term on the right-hand side (RHS) of \eqref{YM-LxF} represents a numerical diffusion and the  flux $\mathcal{F}_{i,j+1}^n$ is obtained using the closure \eqref{YM closure} 
\begin{align}\label{YM-LxF2} 
\mathcal{F}_{i,j+1}^n =   \int_{\Omega} \int_{\mathbb{R}^n} \phi_i(\xi) f(\u)  \mathrm{d}\mu^*_\xi(\u; \vec{u}^n_{j+1}  ) p(\xi) \mathrm{d}\xi,
\end{align} 
where $\vec{u}^n_{j+1} = \left( u_{1,j+1}^n, \dots, u_{N_\xi,j+1}^n \right)^T.$ Furthermore, the time step $\Delta t$ has to be specified. At time $t,$ the CFL condition for equation \eqref{PDE1} is given by  
\begin{equation*}
\max\limits_{ (x,\xi) }\left\{  \sigma\left(  Df\left(u(t,x,\xi) \right) \right) \right\}	\; \Delta t  { = \mbox{\rm{CFL}}}  \Delta x, 
\end{equation*}
where $\sigma$ denotes the spectral radius { of the Jacobian matrix $Df$ and ${\mbox{\rm{CFL}}} \in (0,1].$}  The CFL condition at time $t^n$ needs to be computed in terms of $\vec{u}^{n}_{i,j}$ and the time step $\Delta t$ in the proposed scheme is computed using the following condition:
\begin{align}\label{YM-CFL}
\max\limits_{j=1,\dots,N_x } \left\{ \int_\Omega  \int_{\mathbb{R}^n} \sigma\left(  Df\left(\u \right) \right) \mathrm{d}\mu^*_\xi\left( \u; \vec{u}^n_{j} \right)p(\xi) \mathrm{d}\xi \right\} \Delta t=\mbox{\rm{CFL}}\, \Delta x.
\end{align}

It is a semi-discrete  scheme in the sense that we have only discretized on $(0,T)\times \Bbb R$ without any approximation in  the random space $\Omega.$ 
For a fully discrete scheme,  we need to specify the orthonormal family of basis functions $\{ \phi_i\}_{i=1}^{N_\xi}.$  Furthermore,   in order to evaluate the numerical flux \eqref{YM-LxF2}  suitable numerical quadratures for $\xi \in \Omega$ and  $\u \in \Bbb R^n$ will be applied. The chosen quadratures need  to be also applied to numerically solve the linear programming problem \eqref{ym lin prog} as well as to evaluate the CFL condition \eqref{YM-CFL}. 

\subsection{Summary and Discussion of Proposed Method}
Let us summarize the proposed scheme to compute approximations to the moments $U_i, i=1,\dots,N_\xi$ of the solution $u(t,x,\xi):$ 

For given random initial data $u_0(x,\xi)$ we compute the corresponding Young measure $\nu_{0,x,\xi}$ by means of \eqref{initial-ym-0}. Discrete initial data $u_{i,j}^0$, $i=1,\dots,N_\xi \mbox{ and } j=1,\dots,N_x$ are obtained by equation \eqref{PDE-YM-ID}. Then, the moments $u_{i,j}^n$ are propagated over time for $n=0,1,\dots, N$  according to equation \eqref{YM-LxF} using the flux  \eqref{YM-LxF2}, where $\mu^*_\xi(\vec{u})$ is given by equation \eqref{solution lin prog} and \eqref{ym lin prog}. The time step is chosen at each time step $n$ according to CFL condition \eqref{YM-CFL}. A fully discrete version is presented in \S\ref{sec:colloc} with the corresponding numerical results in \S\ref{sec:num}. Some remarks  are now in order. 

\begin{remark}
The minimization problem \eqref{ym lin prog} is stated for a general parameterized family of measures $\mu_\xi \in \mathcal{P}(\mathbb{R}^n),$ such that the existence of a minimal solution is not necessarily guaranteed.  The results may also depend on the chosen set of basis functions $\{\phi_i \}_{i=1}^{N_\xi}.$  In the next section, we will consider a particular choice of the basis functions $\{\phi_i \}_{i=1}^{N_\xi}$  and discuss expected  properties of the resulting numerical method. As the cost functional the entropy has been chosen and in \S\ref{sec:num} we exemplify that this choice indeed selects the entropy  solution \eqref{PDE2}.

If the flux \eqref{YM closure}, \eqref{solution lin prog} is considered, then the resulting system for the evolution of $\vec{u}$ requires computing the spectral radius of the Jacobian matrix $Df  \in\mathbb{R}^{N_\xi \times n}$ with entries $Df_{i,\ell}$ for  $i=1,\dots,N_\xi$ and $\ell=1,\dots,n$,  
\begin{equation*} 
Df_{i,\ell} = \partial_{u_\ell} \left\{  \int_{\Omega} \int_{\mathbb{R}^n} \phi_i(\xi) f(\u) \; \mathrm{d} \mu^*_\xi(\u; \vec{u} )  p(\xi) \mathrm{d}\xi  \right\}.
\end{equation*}
Unfortunately, in general, it is not expected to find a closed form of $Df_{i,\ell}.$ In the numerical results condition \eqref{YM-CFL} is used.

If we choose particular representations of $\mu_\xi$, some well-known schemes can be recovered on a formal level. For example, if we consider, in \eqref{ym lin prog}, an ansatz for $\mu_\xi$ as $\mu_\xi (\u)= \delta\left(\u - U(t,x,\xi)\right)$ for some unknown function $U=u(t,x,\xi)$, the nonlinear problem \eqref{moment nonlinear} is recovered.  Hence formally, our approach is a generalization of the method \eqref{moment nonlinear} by allowing a more general form of $\mu_\xi$. Alternatively, one could also use the ansatz $\mu_\xi(\u)=\delta\left(\u - \sum\limits_{i=1}^{N_\xi}V_i(t,x) \phi_i(\xi) \right)$ for a set of unknown functions $\{ V_i \}_i$. Then, we  formally recover the classical stochastic Galerkin method, cf.~\cite{chertockchallenges}.  In this case, the degrees of freedom $V_i(t,x)$ are uniquely determined by the set of constraints.

Clearly, the proposed method can be generalized to high-order methods using a suitable polynomial reconstruction and higher-order fluxes in equation \eqref{YM-LxF}. Note that dependence on the cell averages $\vec{u}_j^n$ only enters through the constraints in the linear program to compute the measure $\mu^*_\xi.$  Here, only cell averages have been considered, but clearly, high-order reconstructions of $u$ in spatial direction may also be used here. The same holds true for the higher-order integration in time.  

\end{remark}

\subsection{Fully Discrete Scheme}\label{sec:colloc}
Our aim in this section is to present an example of a fully discrete scheme. Let us assume that $\xi \sim  \mathcal{U}(\Omega)$, i.e.,~$\xi$ is uniformly distributed and  $p(\xi)=p_0 \equiv \frac{1}{|\Omega|}.$ Furthermore, we assume that $\{\xi_i\}_{i=1}^{N_\xi}$ are equidistant points in $\Omega$. To illustrate the strategy for a fully discrete approximation, we will consider piecewise-constant orthonormal basis (\ref{base-coll}), which mimics collocation at equidistant points, i.e.,  
\begin{equation}\label{base-coll}
\phi_i(\xi)=\frac{1}{\sqrt{p_0 \Delta \xi}} \chi_{[\xi_i-\frac{\Delta \xi}2,\xi_i+\frac{\Delta \xi}2 ]}(\xi), \quad i=1,\dots,N_\xi,
\end{equation}  
where $p_0 = \frac{1}{|\Omega|}$ is the constant density of the uniform distribution. This construction yields a set of orthonormal basis functions satisfying
$$
\int_\Omega \phi_i(\xi)\phi_j(\xi)\, p(\xi)\, \mathrm{d}\xi = \delta_{ij}.
$$
The use of piecewise constant, non-overlapping basis functions aligns with the idea of stochastic collocation at equidistant points, while allowing a Galerkin-type discrete formulation. This basis choice simplifies both projection and integration in the stochastic variable, making it particularly suitable for the fully discrete approximation considered here.

We point out that different choices of the random basis functions will yield different fully discrete methods. Furthermore, we consider a discretization with $N_u$ equidistant points $u_\ell \in \Bbb R^n$ of the solution space for $u(t,x,\xi) \in \Bbb R^n.$ According to \eqref{initial-ym-0}, \eqref{PDE-YM-ID} and \eqref{base-coll}, the cell averages of the initial data are  given by 
\begin{align}\label{id-coll}
u_{i,j}^0  =  \frac{1}{\Delta x} \frac{1}{\Delta \xi} \int_{x_j-\frac{\Delta x}{2}}^{x_j+\frac{\Delta x}{2}} \int_{\xi_i-\frac{\Delta \xi}2}^{\xi_i+ \frac{\Delta \xi}2} u_0(x,\xi)  p(\xi)  \mathrm{d}\xi\mathrm{d}x. 
\end{align}
For given moments $u_{i,j}^n$, a discrete approximation $\mu^*_{i,\ell}$ to the measure-valued solution $\mu^*_\xi(\u; \vec{u}^n)$ of the linear programming problem \eqref{ym lin prog} needs to be computed. For a given vector $\vec{u}_j^n = \left( u_{1,j}^n, \dots, u_{N_\xi,j}^n \right)^T$ we approximate 
\begin{equation*}
\mu^*_{\xi_i}(\u_\ell; \vec{u}_j^n) \approx \mu^*_{i,\ell} (\vec{u}_j^n ),
\end{equation*} 
for $i=1,\dots,N_\xi$ and $\ell=1,\dots,N_u.$ The discretization of the linear programming problem \eqref{ym lin prog} reads for the unknowns $\{ \mu_{i,\ell} \}_{i=1,\dots,N_\xi, \ell=1,\dots,N_u}$ 
\begin{subequations} \label{linprog}
\begin{align}
\{ \mu^*_{i,\ell} (\vec{u}_j^n ) \}= & \mbox{ argmin }_{  \{ \mu_{i,\ell} \} }  \Delta \xi \; \Delta u \; 
\sum_{i=1}^{N_\xi}\sum\limits_{\ell=1}^{N_u} \eta(\u_\ell) p_0 \mu_{i,\ell} \label{linprog-aa} \\
\mbox{ subject to }& \mu_{i,\ell} \geq0, \; \forall (i,\ell),  \label{linprog-a} \\
& \mu_{i,\ell} \leq \lambda_F/\Delta u , \; \forall (i,\ell), \label{linprog-c}  \\
& \Delta u \sum_{\ell=1}^{N_u} \mu_{i,\ell} = 1, \; \forall i,  \label{linprog.d} \\
& \Delta u \sum\limits_{\ell=1}^{N_u} \u_\ell \; \mu_{i,\ell} = u_{i,j}^n, \; \forall i. \label{linprog-e}
\end{align}	
\end{subequations}
Here, $0<\lambda_F\le 1 $ is a factor controlling the support of the Young measure and it will be  fixed later in the numerical experiments. The case $\lambda_F=1$ corresponds to the case where atomic solutions are allowed.

The additional constraints guarantee that $\mu^*_{i,\ell}$ is a probability measure for each $i=1,\dots,N_\xi.$ Furthermore, due to the discretization in random space with step size $\Delta \xi,$ there is an upper bound on the size of the measure in each cell. This is also encoded as an additional constraint compared with equation \eqref{ym lin prog}. The stated dependence of $ \mu^*_{i,\ell}$ on $\vec{u}_j^n$ denotes that the solution to the linear program \eqref{linprog} depends only on the current moments on the cell $C_j$ at time $t_n$.  The scheme   is given by equation \eqref{YM-LxF}, where using \eqref{base-coll} and \eqref{linprog}, the flux \eqref{YM-LxF2} reads
\begin{equation} \label{flux-coll}
\mathcal{F}_{i,j}^n =  \Delta u \; \sum\limits_{\ell=1}^{N_u}   f(\u_\ell)  \mu^*_{i,\ell}(\vec{u}_j^n )
\end{equation}
and $\Delta t$ is obtained by the CFL condition \eqref{YM-CFL}
\begin{align}\label{coll-CFL}
\max\limits_{j=1,\dots,N_x } \left\{  \Delta \xi \sum\limits_{i=1}^{N_\xi} \Delta u \sum\limits_{\ell=1}^{N_u} \sigma\left(  Df\left(\u_\ell \right) \right)  \mu^*_{i,\ell} (\vec{u}_j^n )   \; p_0  \right\}\; \Delta t = {\mbox{CFL}}  \Delta x.
\end{align}
\noindent{The equations \eqref{id-coll}, \eqref{linprog}, \eqref{YM-CFL}, \eqref{flux-coll}, and \eqref{coll-CFL}  determine the fully discrete numerical scheme. Before proceeding with the discussion on the scheme properties the following remark on its implementation is in order. }

\begin{remark}
For the solution to the linear programming problem an interior point solver has been used. Therefore, it has been advantageous to replace the equality constraints \eqref{linprog-e} by the following inequalities introducing a possible error $\mathcal{O}(\Delta u)$ in the order of the discretization in $\u.$
\begin{align*}
u_{i,j}^n - \Delta u  \leq  \Delta u \sum\limits_{\ell=1}^{N_u} \u_\ell \; \mu_{i,\ell} \leq  u_{i,j}^n + \Delta u,   \; \forall i. 
\end{align*}	
The dimension of the linear program is $N_u\times N_\xi$ with $2 N_\xi + 2 N_u\times N_\xi$ constraints. The complexity can be reduced to $N_\xi$ linear programming problems of dimension $N_u$ as shown in \S\ref{sec:prop-coll}. While the solution quality  has not been changed, the computational time has been improved. In the numerical simulations presented below the formulation \eqref{linprog2} has been used. 
\end{remark}

\subsection{Properties of the Fully Discrete Scheme}\label{sec:prop-coll}
In this section we will analyze the proposed scheme. Due to the bounds on $\mu_{i,\ell}$  the minimizers exist, if the feasible set is non-empty.  

\begin{lemma}
Provided that the feasible set given by \eqref{linprog-a}--\eqref{linprog-e} is non--empty, the  linear program \eqref{linprog} has a  solution.
\end{lemma}

Due to the particular structure of the constraints and if the entropy $\eta\geq0,$  the problem \eqref{linprog} is equivalent to the following $i=1,\dots,N_\xi$ linear programming problems \eqref{linprog2} in the following sense. 

\begin{lemma}
Assume $\eta(\u)\geq0.$  For each fixed $i=1,\dots,N_\xi$ and a given vector $\vec{u}_j^n $,  the  linear programming problem \eqref{linprog2} has a solution in $\mathbb{R}^{N_u}$ denoted by $\{\nu^*_{i,\ell} \}_{\ell=1}^{N_u}.$ 
\begin{equation} \label{linprog2}
\begin{aligned}
\{ \nu^*_{i,\ell} ( \vec{u}_j^n ) \}= & \mbox{\rm{argmin}}_{  \{ \nu_{\ell} \} }   \Delta u \; \sum\limits_{\ell=1}^{N_u} \eta(\u_\ell) \nu_{\ell} \\
                  \mbox{ \rm{subject to} } & \nu_{\ell} \geq 0, \; \forall \ell, \\
                                      & \nu_{\ell} \leq \lambda_F/\Delta u, \; \forall \ell, \\
                                      & \Delta u \sum_{\ell=1}^{N_u} \nu_{\ell} = 1, \\
		                              & \Delta u \sum\limits_{\ell=1}^{N_u} \u_\ell \; \nu_{\ell} = u_{i,j}^n. 
\end{aligned}	
\end{equation}
Then $\mu^*_{i, \ell} = \nu^*_{i,\ell}$ for all $(i,\ell)$ is a solution to problem \eqref{linprog}. The converse also holds true.
\end{lemma}

Assume that $\lambda_F=1$ and for all  $i = 1, \dots, N_\xi$ it holds that  $\vec{u}^n_{i,j} \in \{ \u_1, \dots, \u_{N_u} \}$, then a solution to the linear programming \eqref{linprog2} is given by
${\nu^*_{i,\ell}} = \frac{1}{\Delta u} \delta_{\ell,\ell_i}$ for $\ell_i$ such that $ \vec{u}^n_{i,j} = \u_{\ell_i}.$  This is the solution expected  for the  stochastic collocation method.  Indeed, assume that $\nu_\ell$ is non-zero for two indices $\ell$ and $k,$ then  due to the fact that $\nu_\ell$ is a probability measure, it holds $\Delta u \nu_k = 1 - \Delta u \nu_\ell.$ Furthermore, due to the equality constraint, it holds $(1 - \Delta u \nu_\ell ) \u_k + \Delta u \nu_\ell = u^n_{i,j}.$ However, since $\eta$ is a convex function, it holds
\begin{equation*}
\eta( \u_{\ell_i} ) = \eta( \vec{u}^n_{i,j}  )  = \eta( (1 - \Delta u \nu_\ell ) \u_k + \Delta u \nu_\ell  \u_\ell ) \leq  \Delta u \eta(\u_k) \nu_k + \Delta u \eta(\u_\ell) \nu_\ell.
\end{equation*}
Therefore, the stochastic collocation solution $\u_{\ell_i}$  is consistent with a solution $\nu^*_{i,\ell}$ to the proposed method.  This proves the following lemma. 
\begin{lemma}\label{lemma uniq}
Assume that $\vec{u}^n_{i,j} \in \{ \u_1, \dots, \u_{N_u} \}$ for each $i.$  Then the discrete measure for all $(i,\ell)$ given by 
\begin{equation*}
\nu^*_{i,\ell}( \vec{u}^n_{j} ) =  \frac{1}{\Delta u} \delta_{\ell,\ell_i} \; \mbox{ for } \ell_i \mbox{ such that } \vec{u}^n_{i,j} = \u_{\ell_i} 
\end{equation*}
is a  solution to problem \eqref{linprog}. If $\eta$ is a strictly convex function, the minimizer is unique.
\end{lemma}
The previous lemma states that we can recover the classical stochastic collocation method under the given assumptions. Note, however, that the method we propose is more general and $\mu^*_\xi(\u; \vec{u}^n_{j+1})$ may be dependent on the full vector of values $\vec{u}^n_{j+1}=(u_{1,j+1}^n,u_{2,j+1}^n,...,u_{N_{\xi},j+1}^n)$ and 
$\mu^*$ may not be a Dirac delta.

Linear programming problems such as \eqref{linprog} and \eqref{linprog2} may have infinitely many solutions. However, in the case of the Burgers equation being one example in \S\ref{sec:num}, the entropy $\eta(u)= \frac{u^2}{2}$ equals the flux $f(u)=\eta(u)$. Since  the value of the cost functional at all minimizers is the same, the numerical flux of the proposed scheme \eqref{YM-LxF} and \eqref{flux-coll}, respectively, is {\em independent} of the particular choice of the minimizer. Under the assumptions of Lemma \ref{lemma uniq} this holds also for general hyperbolic conservation laws. Indeed, Lemma~\ref{lemma uniq} provides the existence of a  unique minimizer $\nu^*_{i,\ell}$ for strictly convex entropies. Hence, for any nonlinear function $\u\to f(\u)$ it holds
\begin{equation*}
f(u^n_{i,j}) = \Delta u  \sum\limits_{\ell=1}^{N_u} f(\u_\ell) \nu^*_{i,\ell}. 
\end{equation*} 
Note that this is a particular setting due to the choice of the basis functions \eqref{base-coll}. In the case of different basis functions or not strictly convex entropies the existence and uniqueness of the minimizer $\mu^*_{i, \ell}$ need to be investigated.

\section{Numerical Results}\label{sec:num}
To illustrate the behavior of the proposed numerical method we present several computational results for the random Burgers equation  as well as for the random  isentropic Euler system. In all simulation results, the random variable is assumed to be uniformly distributed with $\xi \sim \mathcal{U}([-1,1])$ {and} $p(\xi)=1/2.$  We consider one-dimensional spatial domain $[0,1]$ if not stated otherwise. 
 
For the given examples, the conservative variables, the fluxes, the entropies and the CFL condition are specified as follows. In the case of the Burgers equation the nonlinear flux $f(u)=\frac12 u^2$ and the entropy $\eta(u)=\frac12 u^2$ are used in the numerical computation. The time step $\Delta t$ is chosen using the global CFL condition 
\begin{equation}\label{CFL1}
\Delta t = {\rm CFL} \; \frac{  \Delta x}{ \max\limits_{ (x,\xi) } | u_0(x,\xi) | }, 
\end{equation}
where $\rm{CFL}=3/4.$  The initial data $u_0$ are evaluated at the grid points in the $(x,\xi)$ direction.

In the case of the Euler system, the conservative variables are $u=(\rho, q)^T.$ The nonlinear flux is given by $f(u)=(q, \frac{q^2}\rho + \kappa \rho^\gamma)^T.$ We set $\gamma=3/2$ and $\kappa=1$. The entropy for the isentropic Euler system reads   
\begin{equation*}
\eta(u) = \frac12 \frac{q^2}\rho +\frac1{\gamma-1} p(\rho).
\end{equation*}

The eigenvalues of the flux Jacobian matrix are given by 
$\lambda_{\pm}(u) = \frac{q}\rho \pm \sqrt{ p'(\rho) }.$  
Consequently, the CFL condition at discrete points in time $t=n \Delta t$ is 
\begin{equation}\label{CFL2}
\Delta t = {\rm CFL} \; \frac{  \Delta x}{ \max\limits_{ (x,\xi) } \{  |\lambda_+(u(t,x,\xi))|, \; |\lambda_-(u(t,x,\xi)) | \}  }, 
\end{equation}
where ${\rm CFL}=3/4$ and  the maximum is taken over all grid points $(x,\xi)$.

Further, the proposed method requires solving (many) linear programs with equality and inequality constraints. Here, the black-box solver ${\tt linprog}$ of MATLAB R2023b is used. The solver implements two methods and the interior-point method has been used in all examples. The details of the interior-point solver are given in \cite{MR1186163,MR1671584}. The method has been used with standard parameters. All linear programs are solved to optimality.

For a reference solution  of \eqref{PDE1} a standard stochastic collocation method combined with the Lax-Friedrichs finite-volume method in the physical space is used. Let $(\xi_i,x_j)$ be equidistant points for $i=1,\dots,N_\xi$ and $j=1,\dots,N_x$. Let  $C_{i,j}$ denote grid cells of size $\Delta \xi \times \Delta x$ with centers at $(\xi_i,x_j)$. Denote the cell averages by  
\begin{equation*} 
v^0_{i,j} = \frac{1}{ \Delta x \; \Delta \xi } \int_{C_{i,j}} u_0(x,\xi) p(\xi) \mathrm{d}\xi\mathrm{d}x , \; i=1,\dots, N_\xi, \; j=1,\dots,N_x. 
\end{equation*}
The stochastic collocation Lax-Friedrichs finite-volume approximation  to the cell averages of the solution $u(t,x,\xi)$ to \eqref{PDE1}  at time $t_n=n \Delta t$ is given by 
\begin{equation} \label{coll}
v^{n+1}_{i,j} = \frac12  \left(  v^n_{i,j+1} + v^n_{i,j-1} \right) - \frac{\Delta t}{2 \Delta x} p(\xi_i)  \left(  f\left( \frac{v^n_{i,j+1} }{p(\xi_i)} \right)  -  f\left( \frac{v^n_{i,j-1}}{p(\xi_i)} \right) \right).
\end{equation}
High-order stochastic collocation approximations are also possible and we refer to \cite{MR4772617} for more details and further examples. The scheme \eqref{coll} is the simplest first-order scheme presented therein. The time step $\Delta t$ is chosen according to the CFL condition \eqref{CFL1} or \eqref{CFL2}, respectively.

Unless stated otherwise, all comparisons of numerical solutions are done with respect to the reference solution  $v^n_{i,j}$ obtained by the stochastic collocation Lax-Friedrichs finite-volume method \eqref{coll} on the same grid in $(x,\xi)$ and the same CFL condition \eqref{CFL1} or \eqref{CFL2}. The errors are computed in the $L^1$-norm  for a given terminal time  $T$ stated explicitly in equation \eqref{L1err}. 

\subsection{Convergence Studies for the Burgers Equation}
The spatial computational domain is set to $[0,1], $ $\Omega = [-1,1] $ and the phase space is restricted to $\u \in [-5,5]$. The latter is discretized by $N_u = 100$ equidistant cells.

Initial data are chosen as follows   
\begin{equation}\label{id-figure1}
u_0(x,\xi) = \chi_{[0, 1/2]}(x) \; \xi.   
\end{equation}
%The reference solution $v^n_{i,j}$  is computed by  the stochastic collocation Lax-Friedrich finite volume method \eqref{coll}. 
The time step is chosen according to the CFL condition \eqref{CFL1}.  In all cases, we compute the error  between the numerical solution $u^n_{i,j}$ and the reference solution  $v^n_{i,j}$ at time $T=1/2,$ i.e., 
\begin{equation}\label{L1err}
err = \Delta x \Delta \xi \sum\limits_{i,j} | u^N_{i,j} - v^N_{i,j} |,
\end{equation}
where $N\Delta t=T.$   In Table~\ref{tab1}, we present the experimental order of convergence for the approximation in physical space. Consequently, the grid in 
$(u,\xi)$ has been fixed and the spatial grid has been successively refined. As expected, the first-order convergence rate is observed. In Table~\ref{tab2} the convergence  for the approximation in the random direction is presented; thus the grid in the $(x, \u)$-direction is fixed. The observed experimental error  is of order $\mathcal{O}(10^{-3})$. This is the error of the Lax-Friedrichs finite-volume method  with  a spatial discretization  of $N_x=500$ points on $[0,1].$

\begin{table}[ht!]
\centering
\begin{tabular}{ccc}
\toprule
$N_x$ & $L^1([0,1]\times [-1,1])$ & Rate \\
\hline
40 & 1.3697E-02 & 8.0190E-01 \\  
60 & 9.7416E-03 & 8.4042E-01 \\  
80 & 7.5849E-03 & 8.6986E-01 \\  
100 & 6.2061E-03 & 8.9912E-01 \\  
120 & 5.2575E-03 & 9.0976E-01 \\  
140 & 4.5704E-03 & 9.0857E-01 \\  
160 & 4.0381E-03 & 9.2726E-01 \\  
180 & 3.6182E-03 & 9.3241E-01 \\  
200 & 3.2794E-03 & 9.3315E-01 \\  
\bottomrule 
\end{tabular}
\caption{\sf Spatial grid with  $N_x$ equidistant points distributed on the physical domain $[0,1].$ Reported is the $L^1$-norm of the error with respect to $(x,\xi)$ with initial data \eqref{id-figure1} and simulated up to terminal time $T=1/2.$ In $\xi$-direction  $N_\xi=5$ and in $\u$-direction $N_u=10$ discretization points have been used.}
\label{tab1}
\end{table}

\vskip 20pt 
\begin{table}[ht!]
		\centering
		\begin{tabular}{ccc}
			\toprule
			$N_\xi$ & $L^1([0,1]\times [-1,1])$ & Rate \\
			\hline
			40  & 1.1493E-02 & 1.0962 \\  
			60  & 7.5345E-03 & 1.0414 \\  
			80  & 5.6034E-03 & 1.0294 \\  
			100 & 4.4599E-03 & 1.0228 \\  
			120 & 3.7040E-03 & 1.0187 \\  
			140 & 3.1671E-03 & 1.0158 \\  
			160 & 2.7662E-03 & 1.0137 \\  
			180 & 2.4553E-03 & 1.0121 \\  
			200 & 2.2073E-03 & 1.0108 \\  
			\bottomrule 
		\end{tabular}
\caption{\sf Spatial grid with  $N_\xi$ equidistant points distributed on the physical  domain $[-1,1]$. Reported is the $L^1$-norm of the error with respect to $(x,\xi)$ with initial data \eqref{id-figure1} and simulated up to terminal time $T=1/2.$ In $x$-direction  $N_x=500$ and in $\u$-direction $N_u=10$ discretization points have been used. }\label{tab2}
\end{table}

\subsection{Burgers Equation with Sinusoidal Initial Data}

\begin{figure}[ht!]
\centerline{\includegraphics[trim=0.cm 0.6cm 1.cm 0.7cm, clip, width=7.cm]{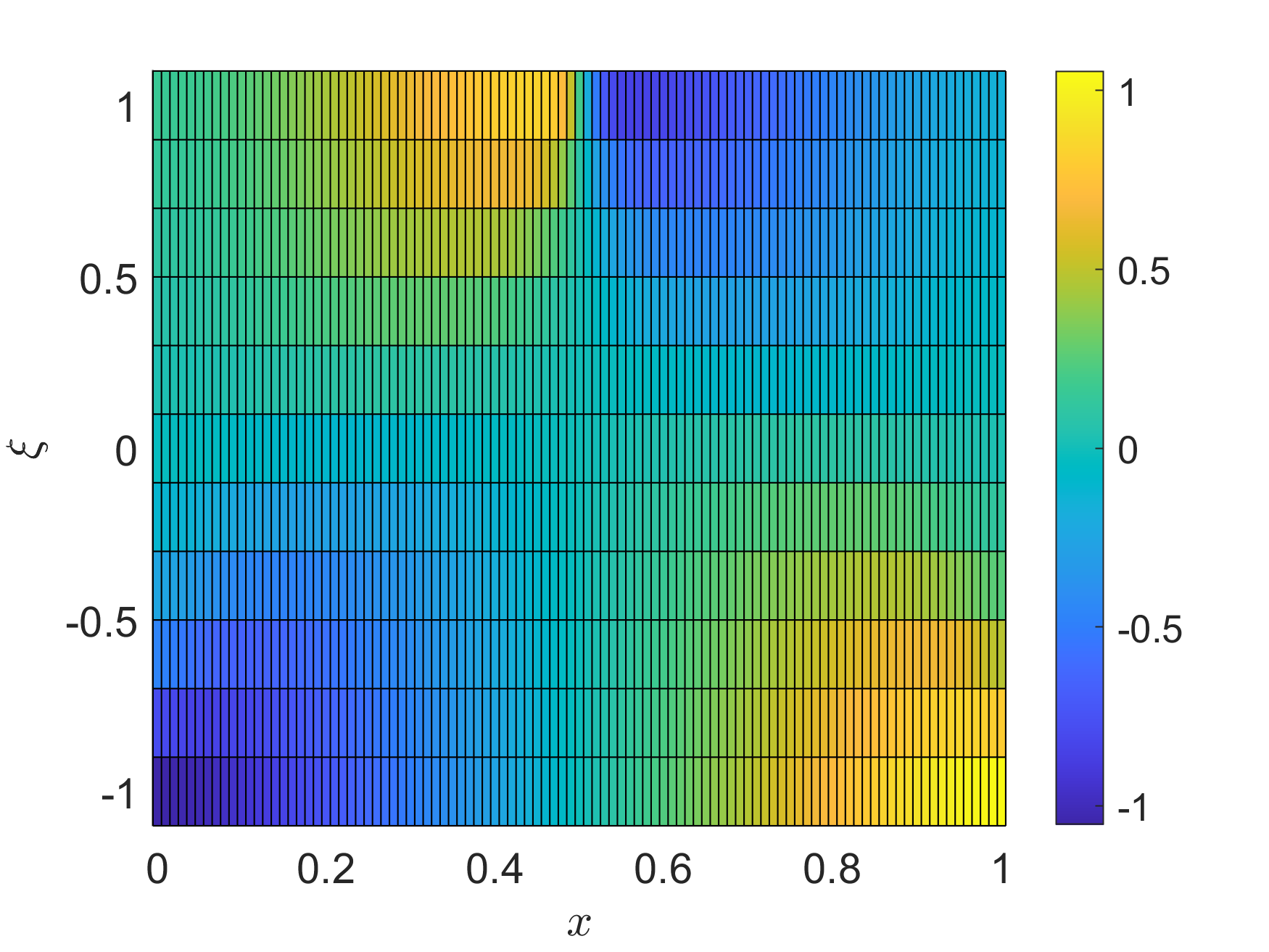} \hspace{1cm}
            \includegraphics[trim=0.cm 0.6cm 1.cm 0.7cm, clip, width=7.cm]{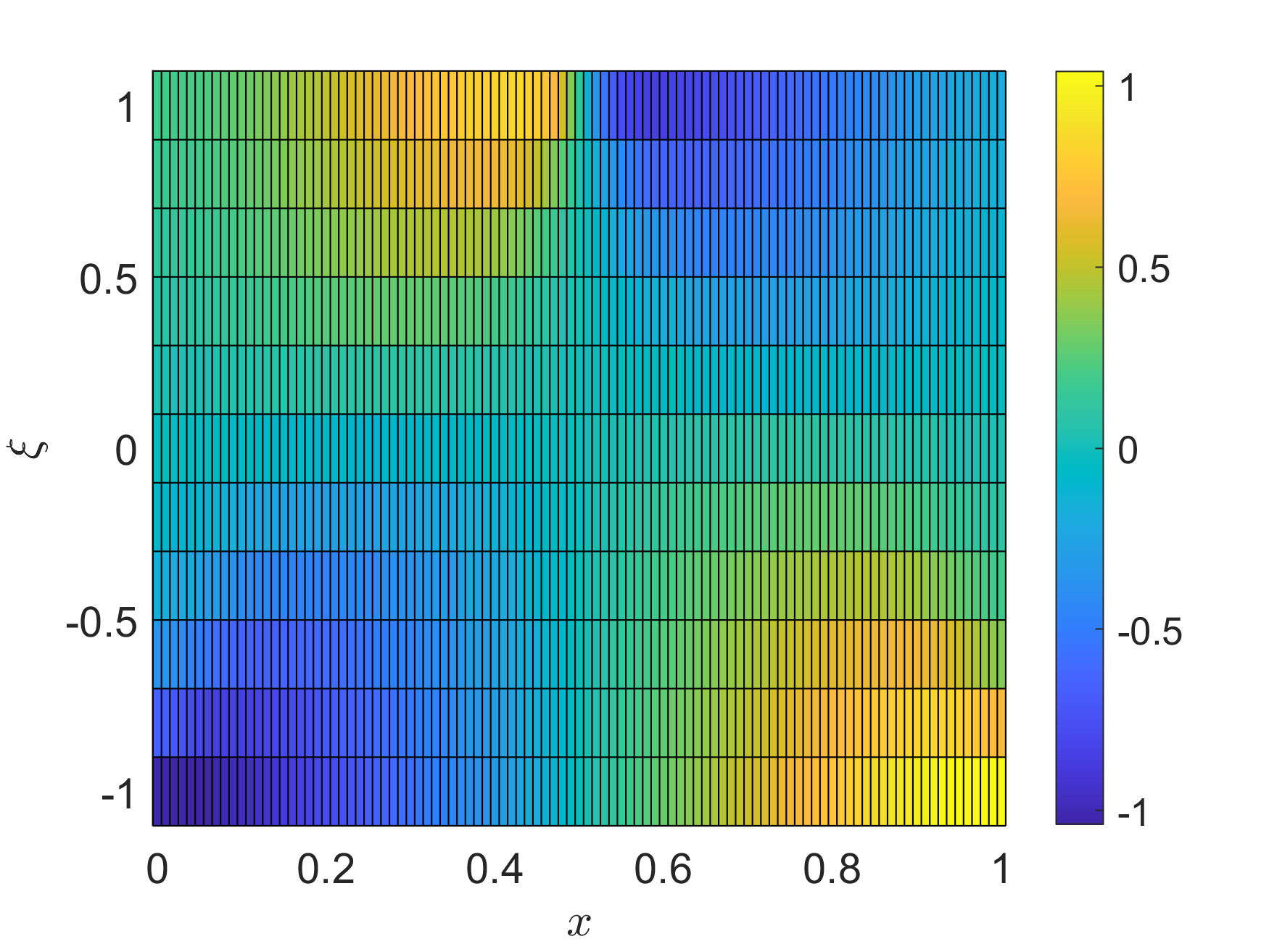}}
\caption{\sf Numerical solutions obtained by  the stochastic collocation Lax-Friedrichs method (left) and  the proposed method (right). Initial data are $u(0,x,\xi)=\xi \sin(2 \pi x).$ The final time $T=1/4.$}
\label{fig1}
\end{figure}

In this experiment we consider the initial data 
\begin{equation*}
u_0(x,\xi) = \xi \sin\left(2\pi x\right).
\end{equation*}
The final time is set to   $T=1/4.$   Although the initial data are smooth, a shock will be formed at a  later time $T>0$ for any  $\xi>0.$  For $\xi<0$, the solution remains smooth. The spatial grid has $N_x=100$ points and the grid in the random space has $N_\xi=10$ points. For the discretization in phase space $\u \in [-5,5]$,  $N_u=100$ equidistant points are used. The time step $\Delta t$ is the same for both schemes. The error in the $L^1$-norm in $(x,\xi)$ at the final time between the reference solution $v_{i,j}^N$ and the numerical solution $u_{i,j}^N$  is $err=4.9e$-$04$, where $err$ is given by equation \eqref{L1err}.  A visualization of the numerical solution obtained by both schemes is given in Figure \ref{fig1}. While both solutions seem similar,  we  observe that the proposed method is slightly more diffusive compared to the stochastic collocation Lax-Friedrichs finite-volume method. This is visible for $\xi$ close to one at the center of the domain $x=1/2.$ 

\begin{figure}[ht!]
\centerline{\includegraphics[trim=0.cm 0.6cm 1.cm 0.7cm, clip, width=7.cm]{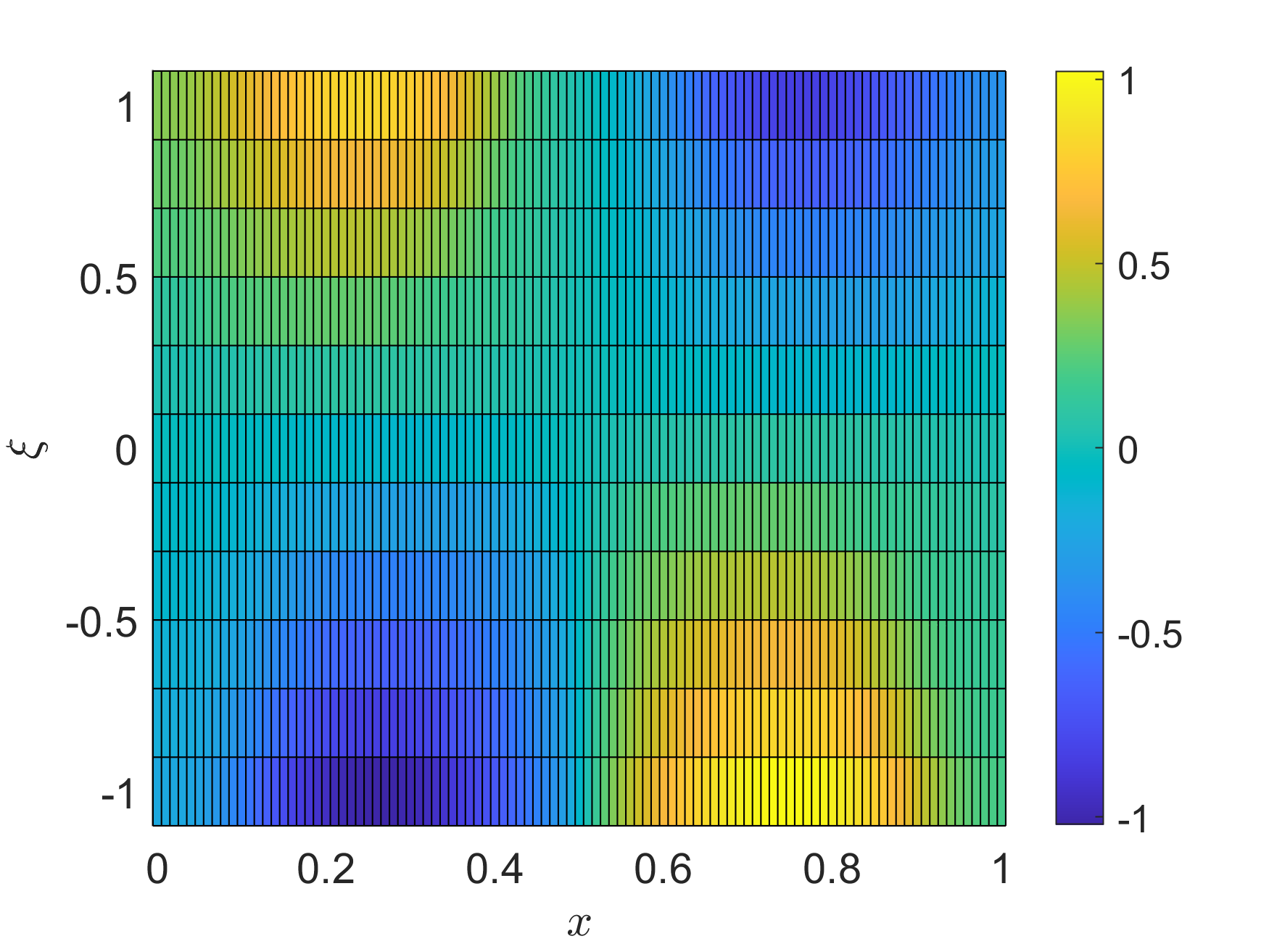}}
\caption{\sf Numerical solution obtained by the proposed method with modified cost functional in the linear program \eqref{linprog}. Initial data are  $u(0,x,\xi)=\xi \sin(2 \pi x).$ The final time is $T=1/4.$}
\label{fig2}
\end{figure}

In the presented method, the entropy is minimized to obtain an approximation to the Young measure. This is necessary to obtain the entropy solution as the following example illustrates. If, for example, we change the minimization problem to feasibility only, i.e., we set $\eta(\u)=1$ in the cost functional of problem  (\ref{linprog}),  the result in Figure \ref{fig2} is obtained. Comparing with  Figure \ref{fig1} we observe that in this case the scheme selects the wrong weak solution consisting of  rarefaction waves even in the case $\xi>0.$ This is not the correct weak entropy solution. 

\subsection{Burgers Equation with Non-Atomic Support}
In this subsection, we consider an example given for the (deterministic) Burgers equation. This example is inspired by \cite{Tadmor-1Dburgers}. In this example, the Young measure may have non-atomic support but its first-order moment solves a Riemann problem with a shock.  This solution can be also obtained by our framework when considering absolutely continuous measures with respect to the Lebesgue measure instead of the whole probability space $\mathcal{P}(\mathbb{R}^n)$. We set the factor $\lambda_F$ in \eqref{linprog-c} to be 0.05. In all runs of this subsection, we do not prescribe any parametric Young measure at $t=0$; we initialize only with the first moment via (\ref{2.14aaa}), and the Young measure is then generated by the first time step of the linear program under the chosen $\lambda_F$. By modifying the upper bound on the value of the weights of the Young measure, the atomic solution is excluded and the linear programming solver requires to find a different Young measure. 
With this small modification, we are able to obtain a solution similar to the one given in  \cite{Tadmor-1Dburgers}. 

\begin{figure}[ht!]
\centerline{\includegraphics[trim=0.2cm 0.2cm 1.1cm 0.7cm, clip, width=7.cm]{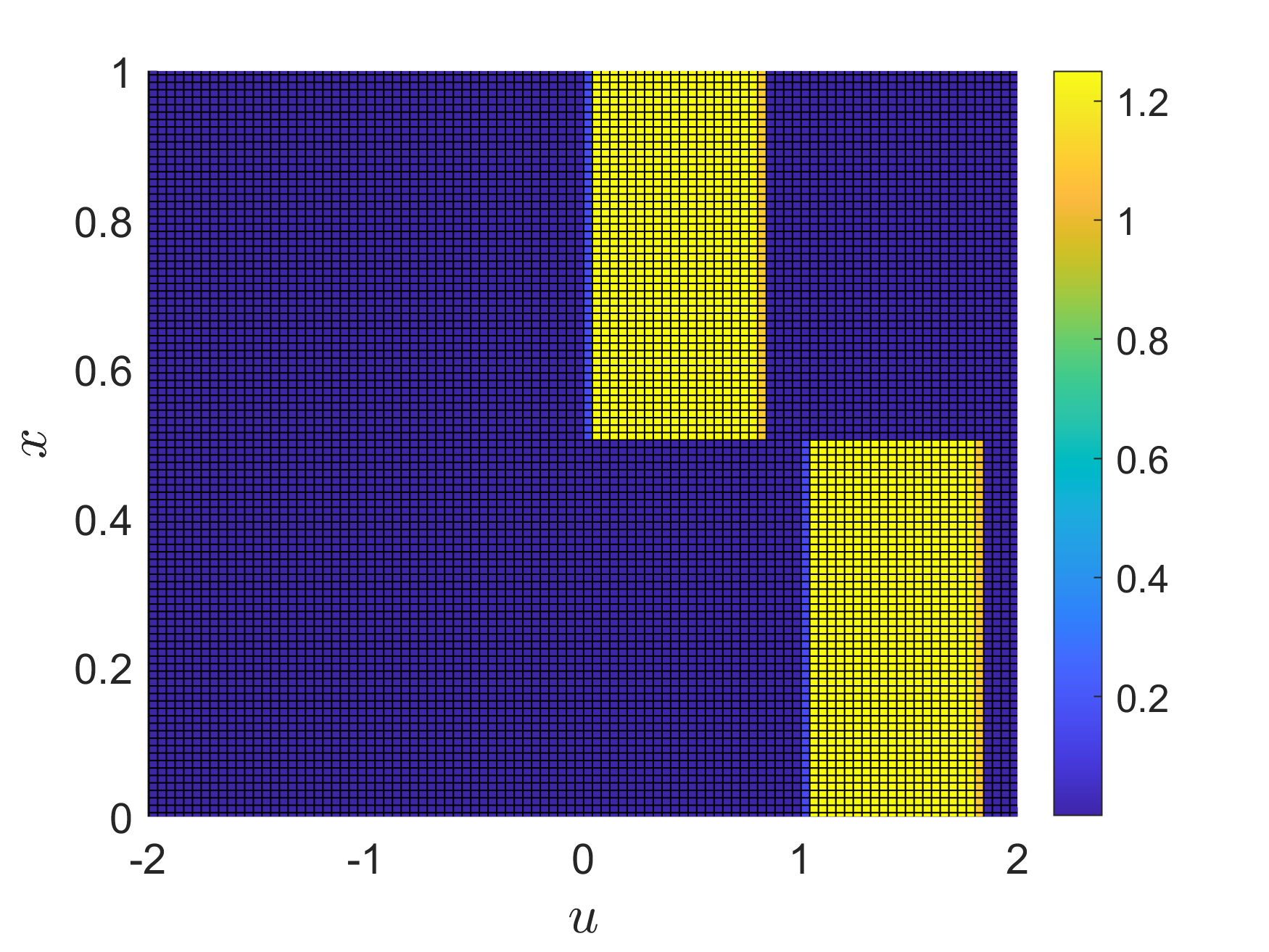} \hspace{1cm}
            \includegraphics[trim=0.2cm 0.2cm 1.1cm 0.7cm, clip, width=7.cm]{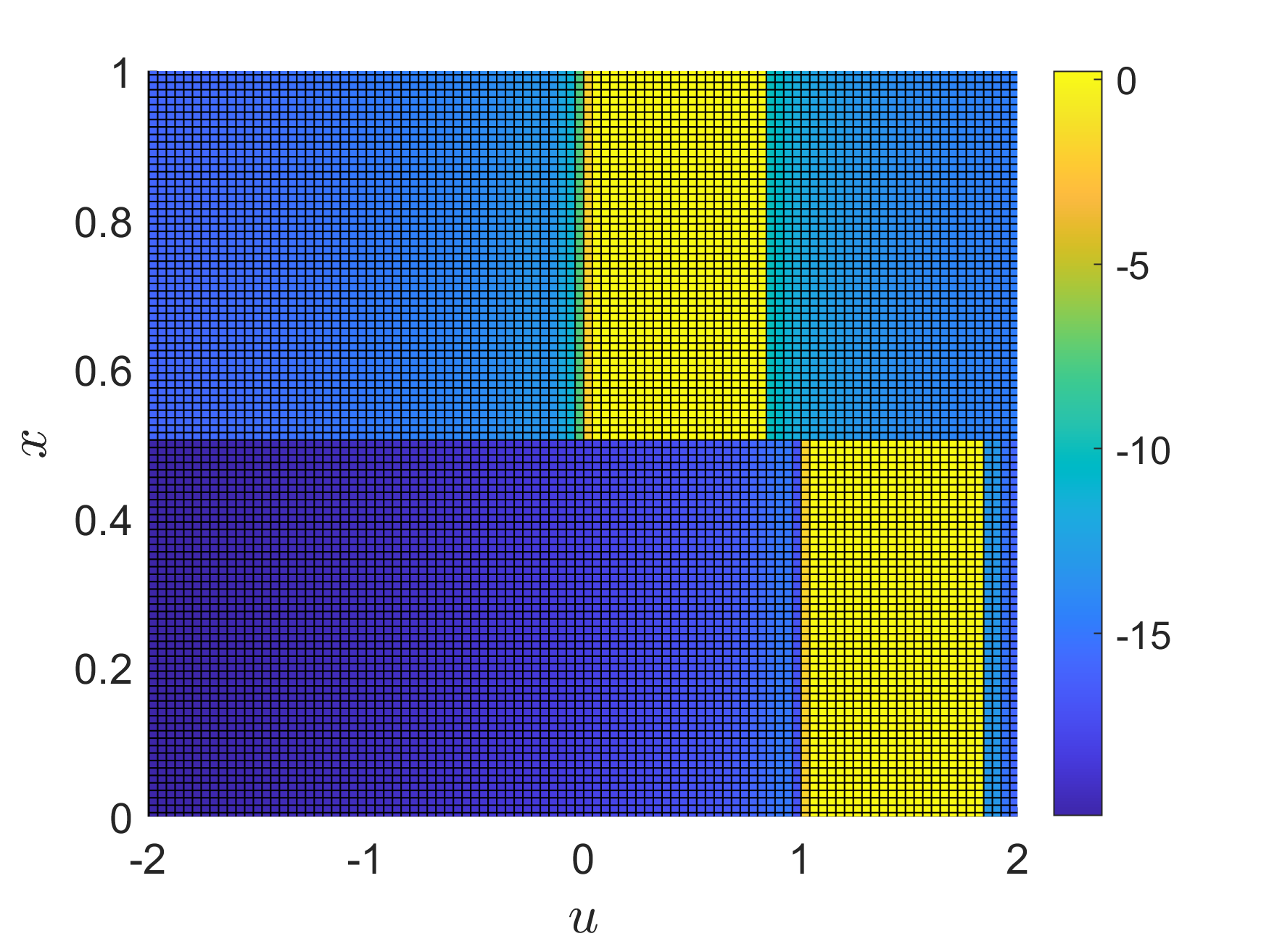}}
\caption{\sf  Young measure at the initial time $t = 0$ (left) and its logarithmic values (right).\label{fig111}}
\end{figure}

In this example, we consider the initial data 
$$
u_0(x)=\begin{cases}
             1.5 & \mbox{if } x< 0.5, \\
             0.5 & \mbox{otherwise},
           \end{cases}
$$
\noindent{subject} to the free boundary conditions imposed in the computational domain $x\in [0,1]$. In Figure \ref{fig111}, we plot the support of the Young measures at the initial time. 

\begin{figure}[ht!]
\centerline{\includegraphics[trim=0.2cm 0.2cm 1.1cm 0.7cm, clip, width=7.cm]{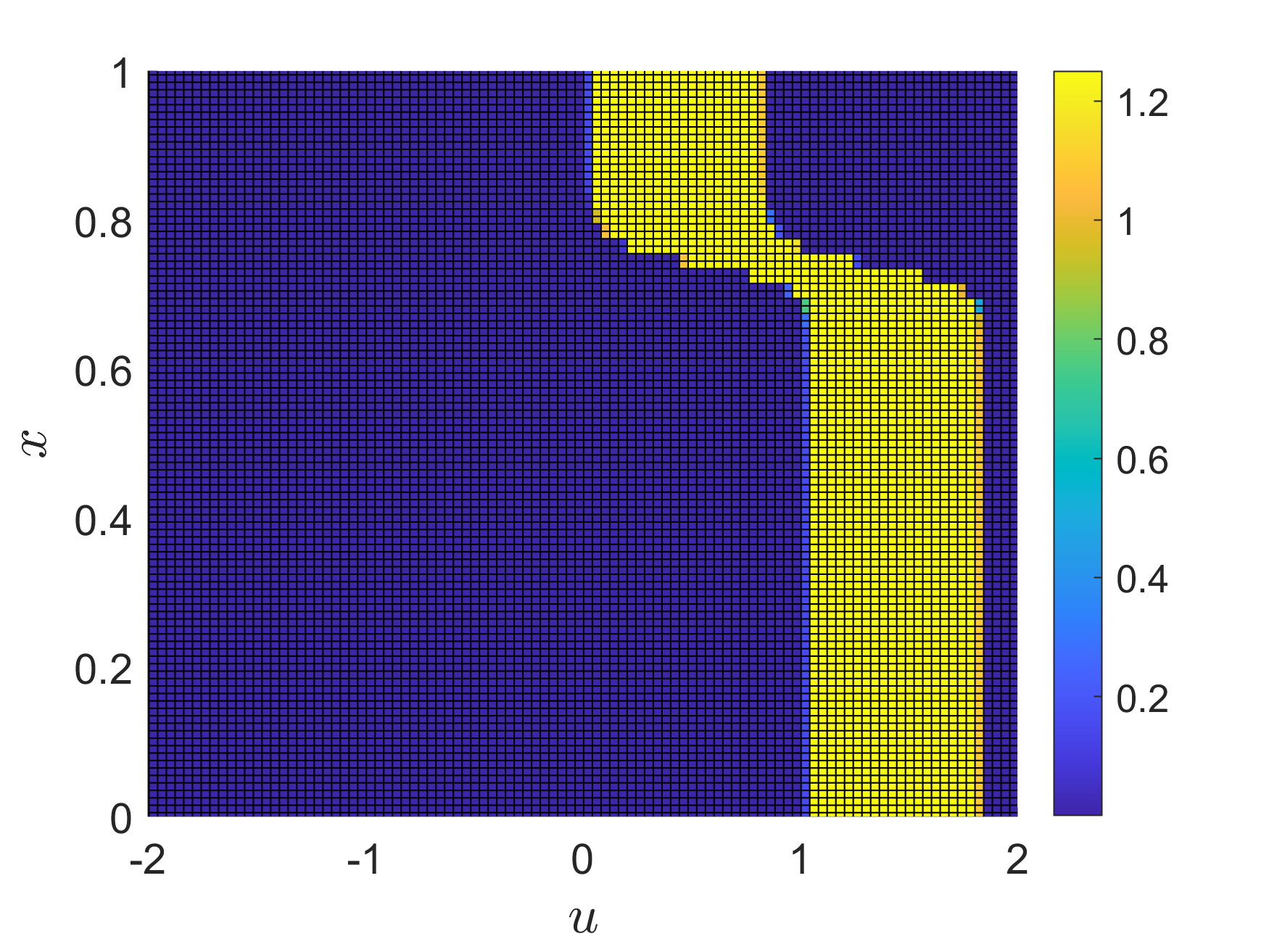} \hspace{1cm}
            \includegraphics[trim=0.2cm 0.2cm 1.1cm 0.7cm, clip, width=7.cm]{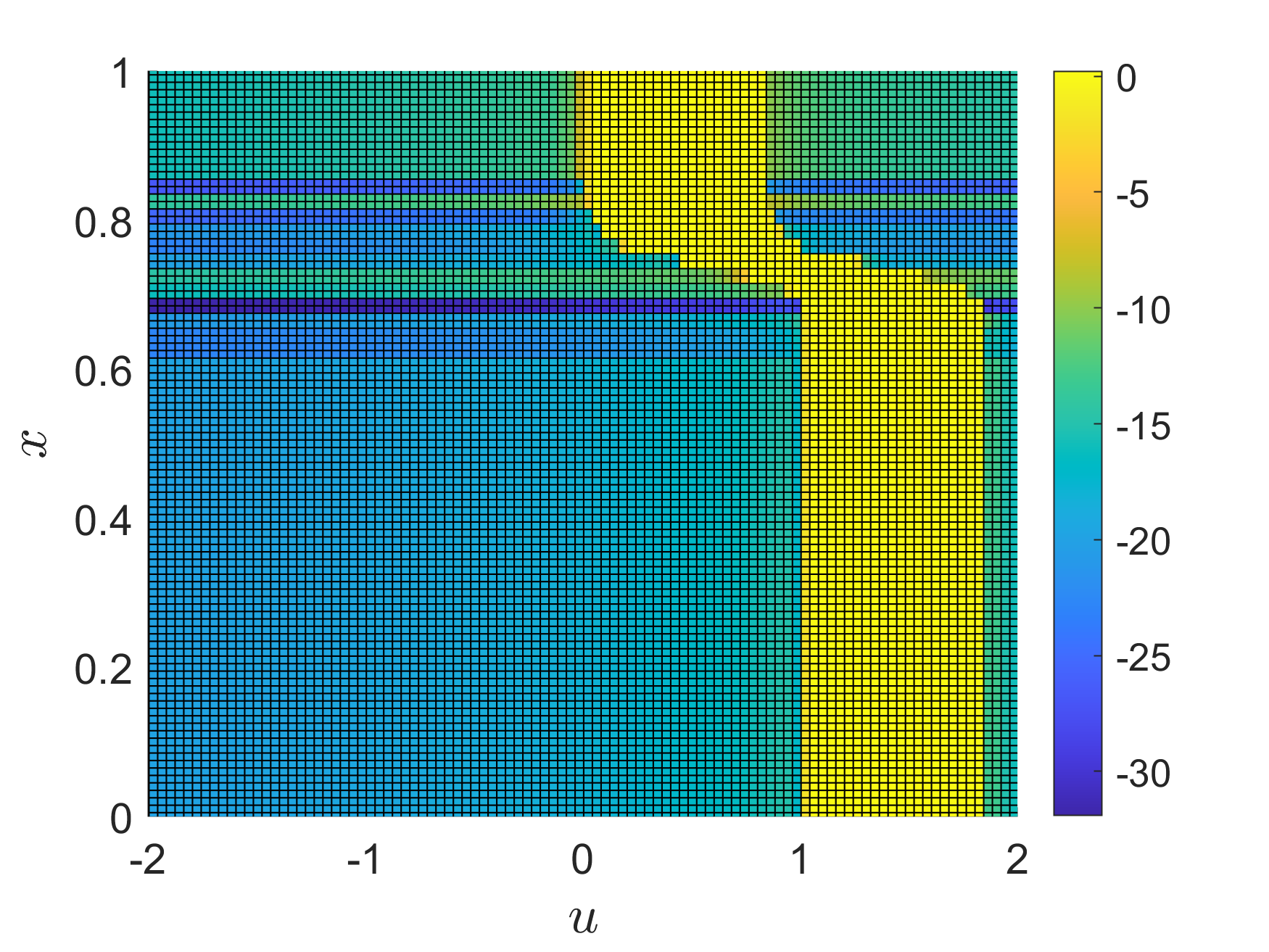}}
\caption{\sf Young measure at the final time $T = 0.25$ (left) and its logarithmic values (right) with $\lambda_F=0.05$.\label{fig222}}
\end{figure}

We compute the numerical results up to the final time $T = 0.25$ with the spatial domain discretized using $N_x = 100$ grid points. The phase space is confined to the interval $\mathbf{u} \in [-2, 2]$ and is uniformly discretized into $N_u = 100$ equidistant cells. The obtained numerical results are presented in Figure \ref{fig222}. From this example, we observe that the measure with non-atomic support can be captured by the linear programming algorithm if we choose $\lambda_F\neq 1$.

\begin{figure}[ht!]
\centerline{\includegraphics[trim=0.2cm 0.2cm 1.1cm 0.7cm, clip, width=7.cm]{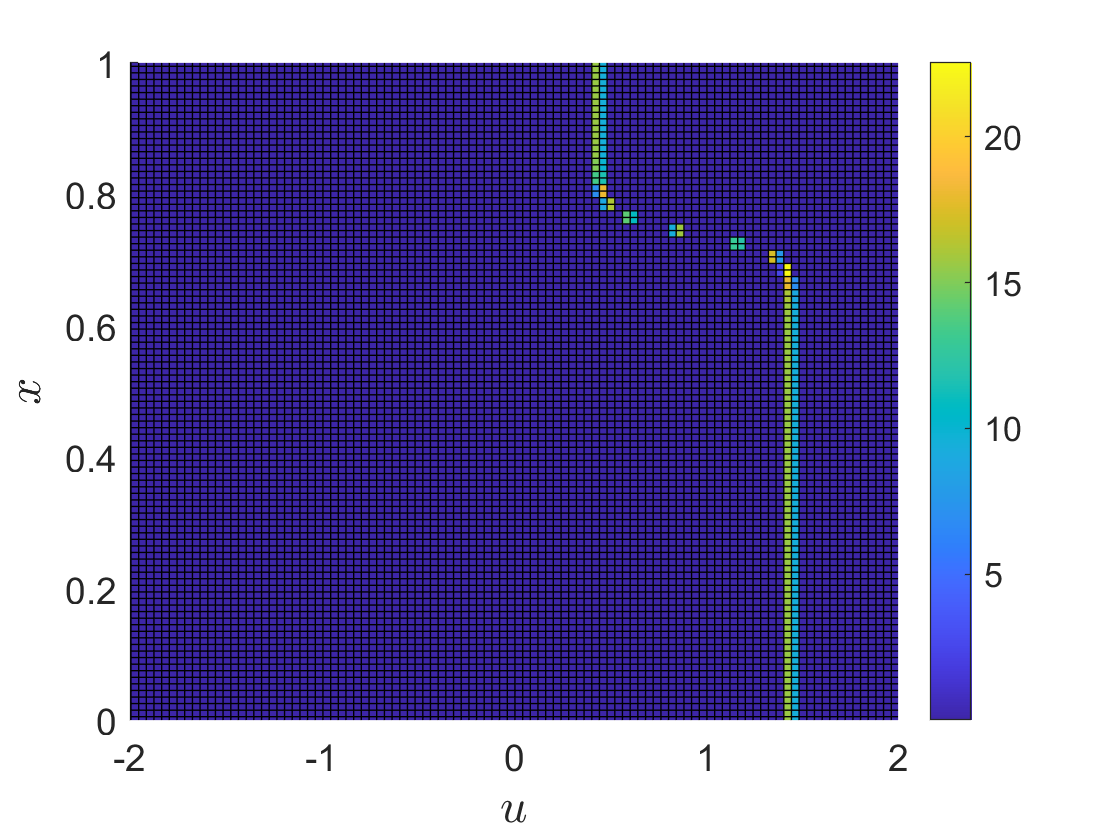} \hspace{1cm}
            \includegraphics[trim=0.2cm 0.2cm 1.1cm 0.7cm, clip, width=7.cm]{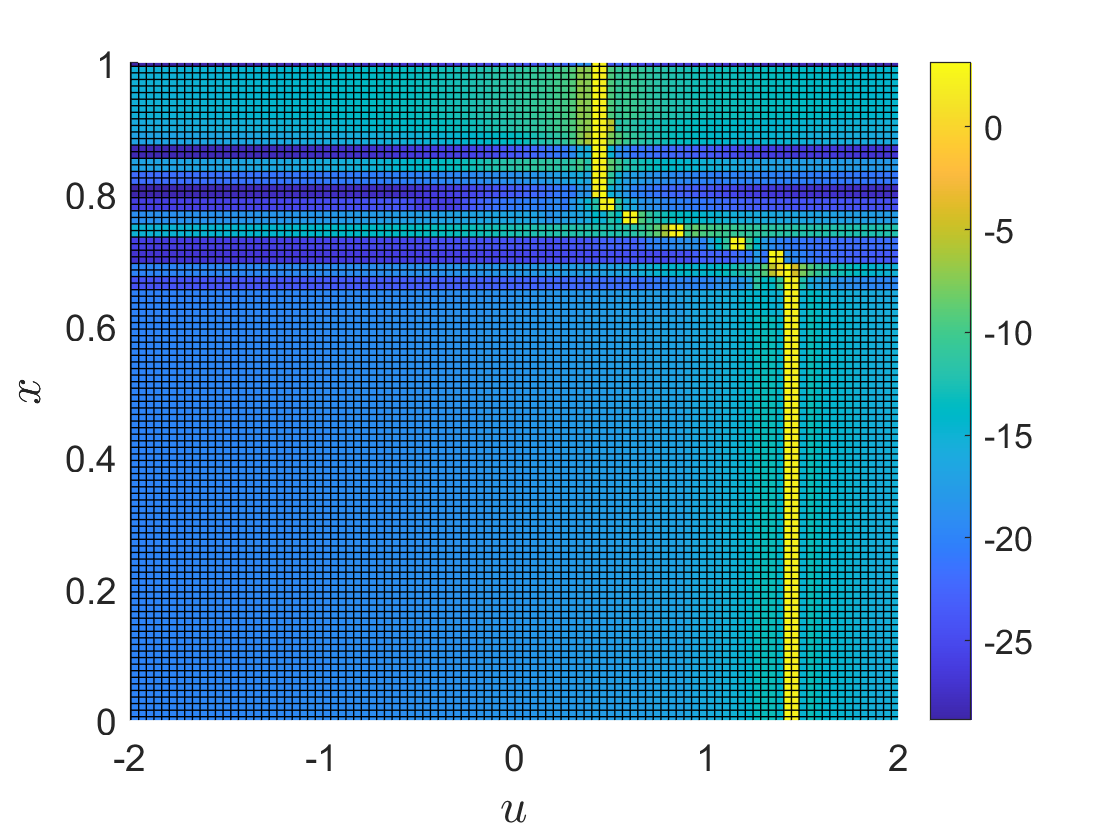}}
\caption{\sf  Young measure at the final time $T = 0.25$ (left) and its logarithmic values (right) with $\lambda_F=1$.}\label{fig222a}
\end{figure}

We also show the computed numerical result with $\lambda_F=1$ in Figure \ref{fig222a}, where one can see the differences between the solutions when using different values of $\lambda_F$. As expected, for $\lambda_F=1$, minimization occurs over the full feasible set and the solution becomes atomic; see Figure \ref{fig222a} for a side-by-side contrast with the non-atomic case $\lambda_F=0.05$.

\subsection{Riemann Problem for the Isentropic Euler System}

The proposed method is  applicable also for systems of hyperbolic conservation laws. We now consider a random isentropic Euler system and take an equidistant grid in $(x,\xi) \in [-1,1]\times[-1,1]$ with $N_x=100$ and $N_\xi=10$ points. 

\begin{figure}[ht!]
\centerline{\includegraphics[trim=0.7cm 0.6cm 1.2cm 0.7cm, clip, width=7.cm]{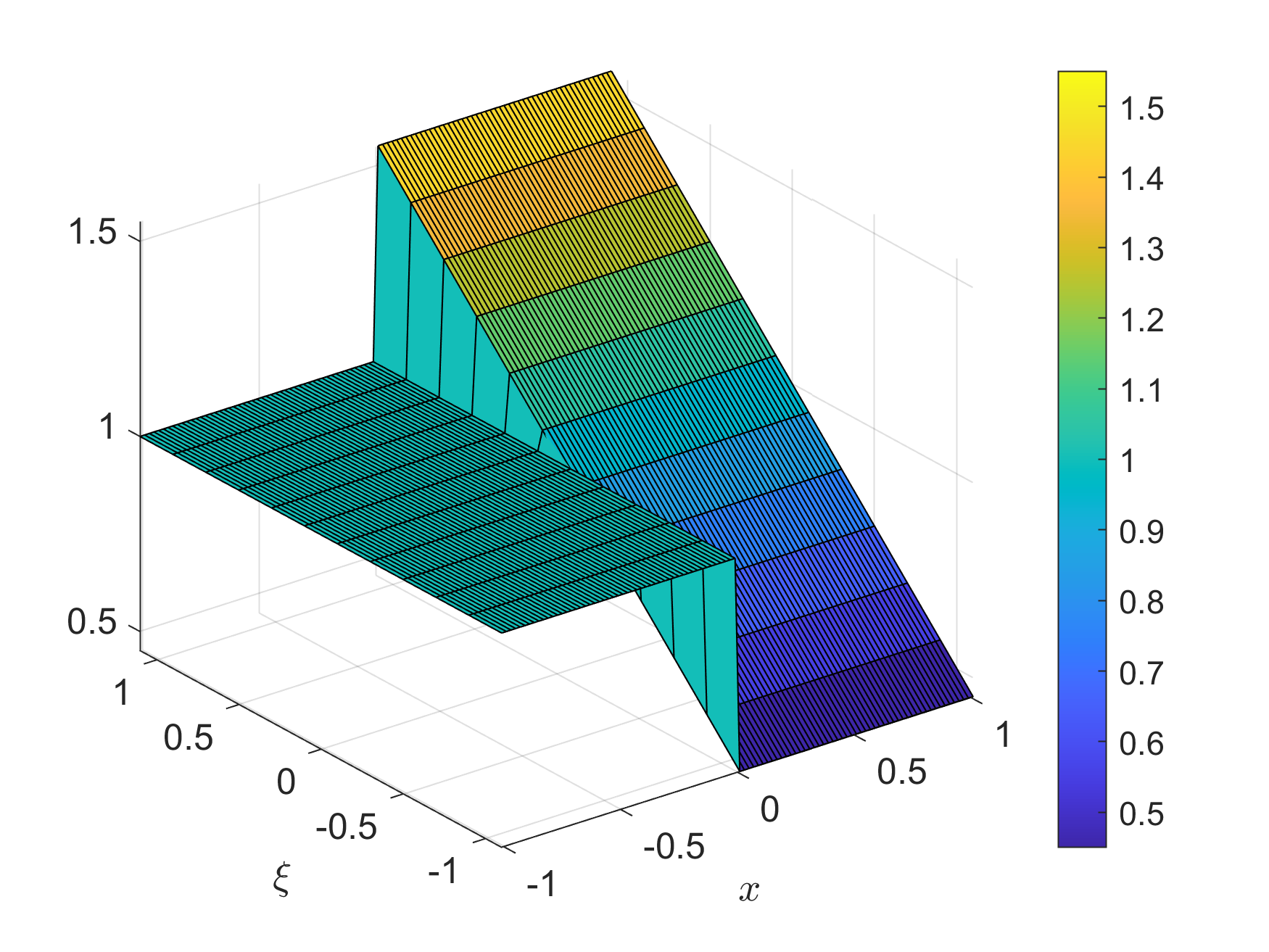} \hspace{1cm}
            \includegraphics[trim=0.7cm 0.6cm 1.2cm 0.7cm, clip, width=7.cm]{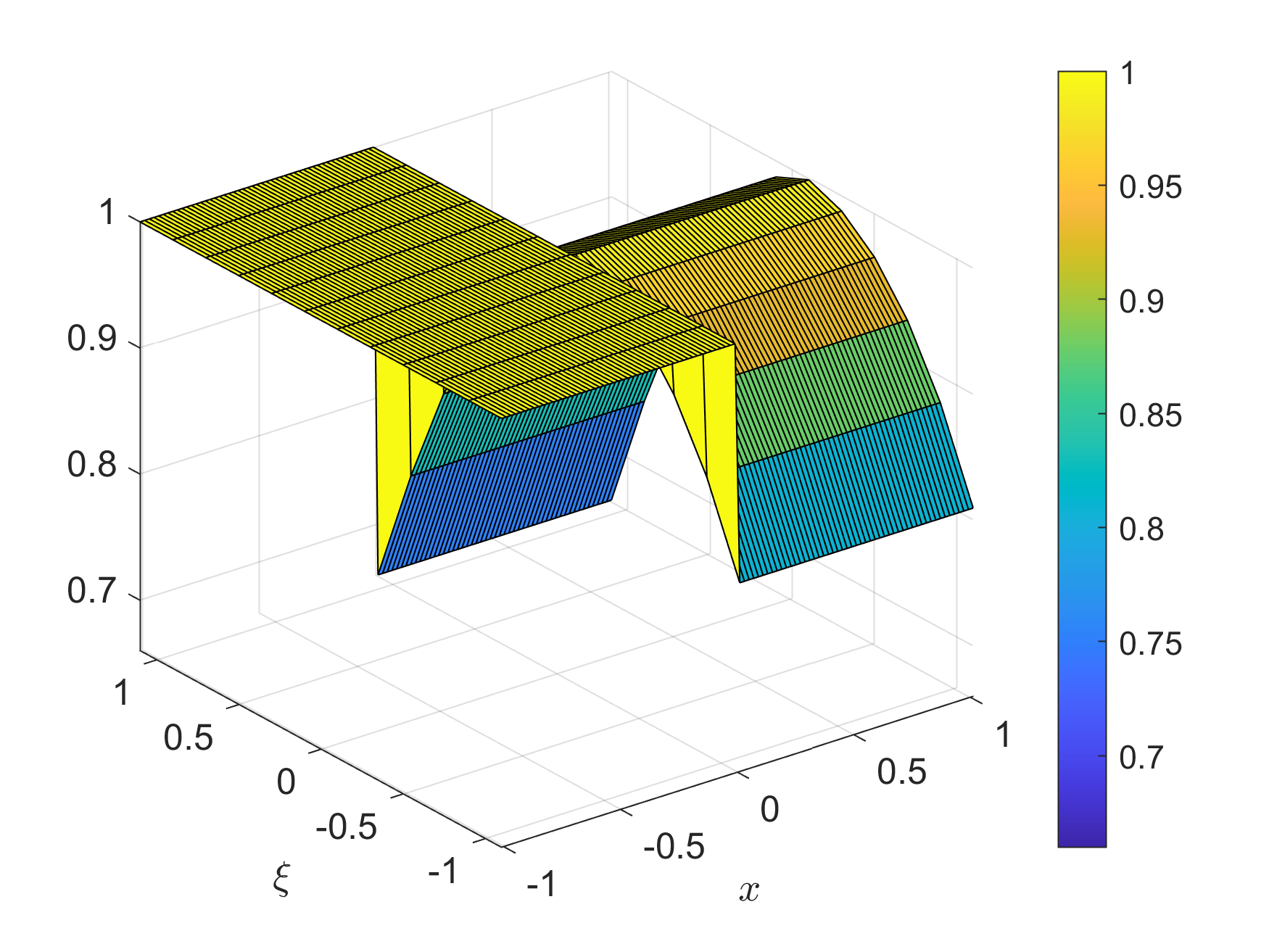}}
\vskip 20pt
\centerline{\includegraphics[trim=0.7cm 0.6cm 1.2cm 0.7cm, clip, width=7.cm]{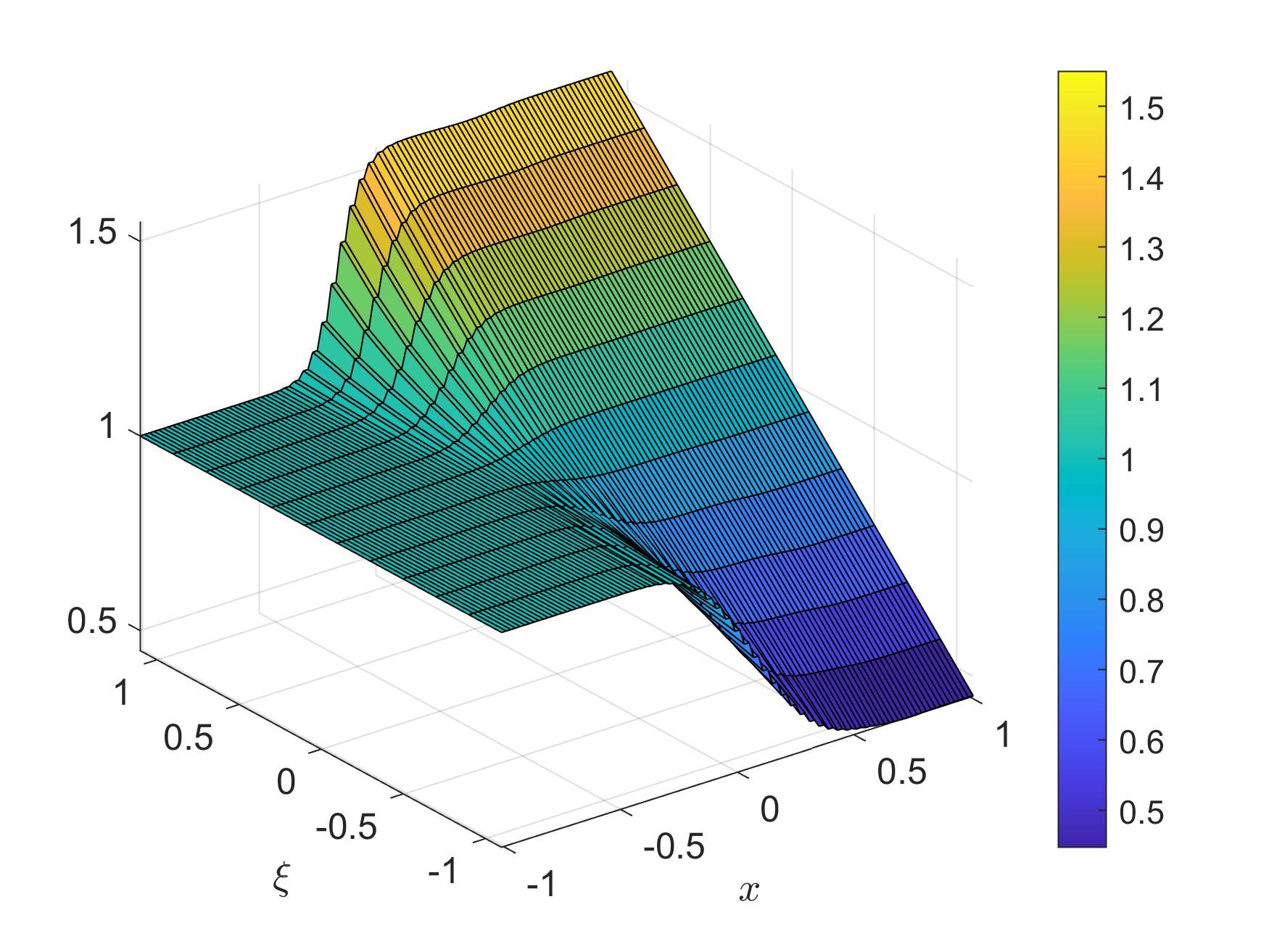} \hspace{1cm}
            \includegraphics[trim=0.7cm 0.6cm 1.2cm 0.7cm, clip, width=7.cm]{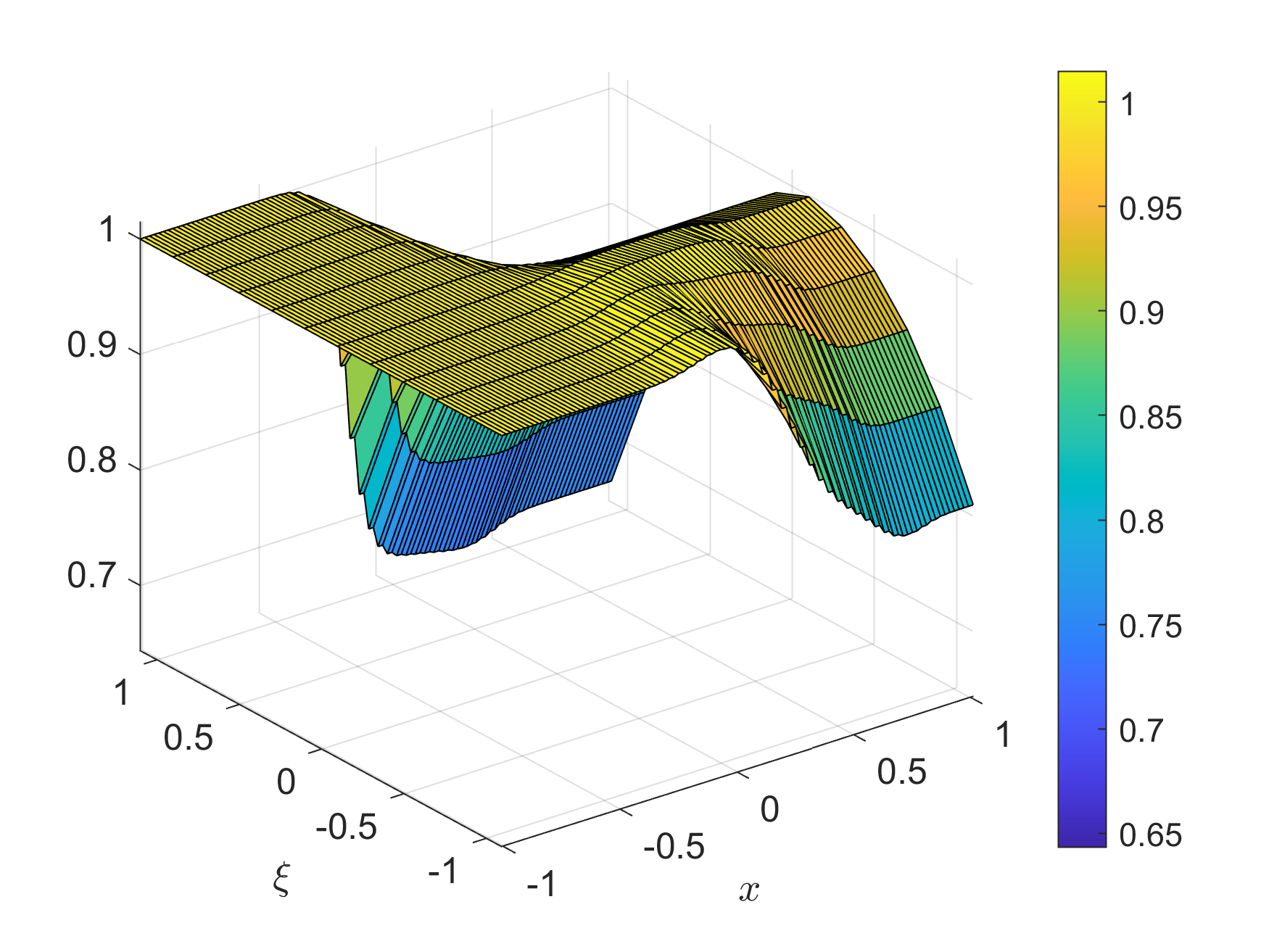}}
\vskip 20pt
\centerline{\includegraphics[trim=0.7cm 0.6cm 1.2cm 0.7cm, clip, width=7.cm]{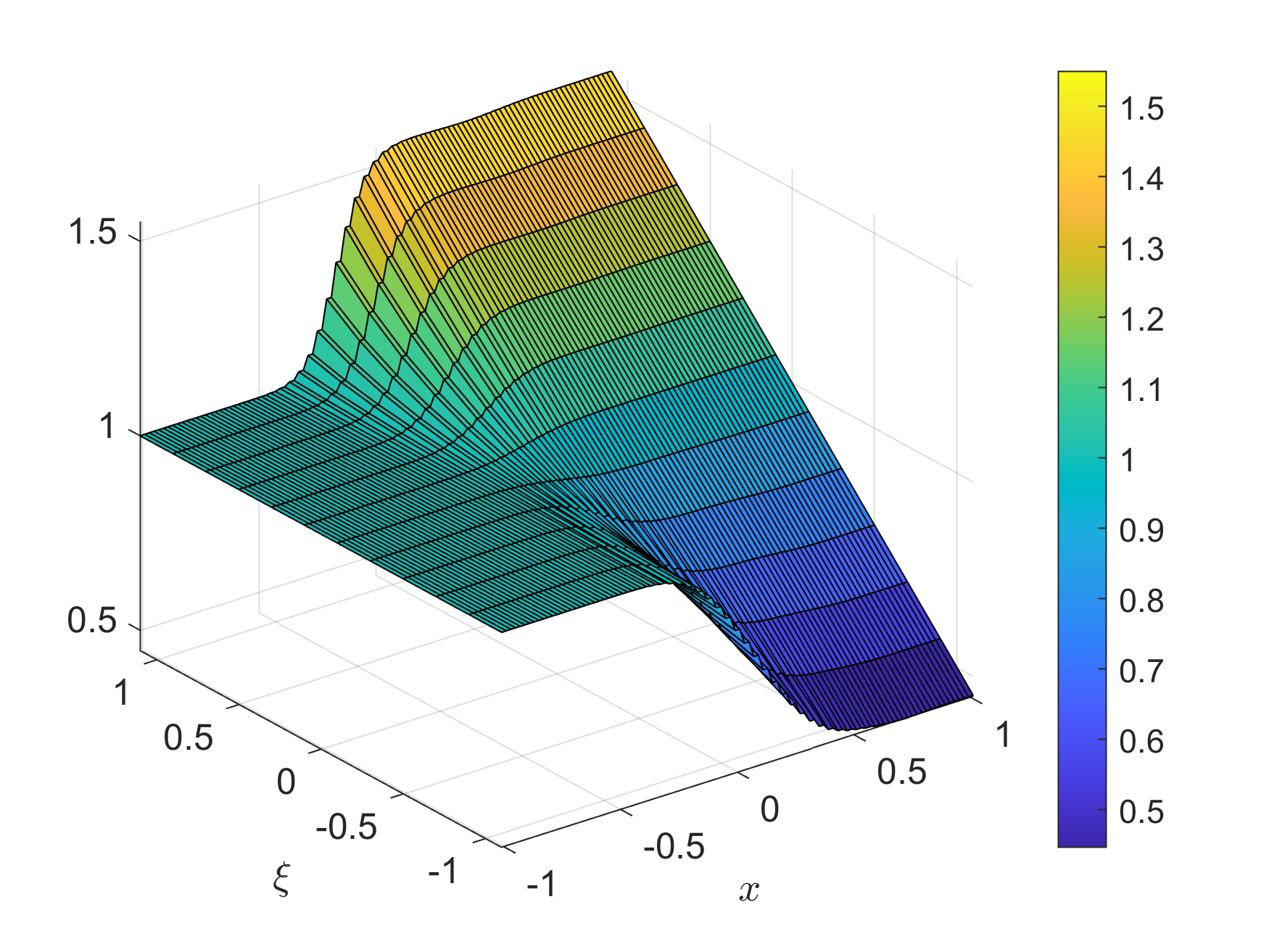} \hspace{1cm}
            \includegraphics[trim=0.7cm 0.6cm 1.2cm 0.7cm, clip, width=7.cm]{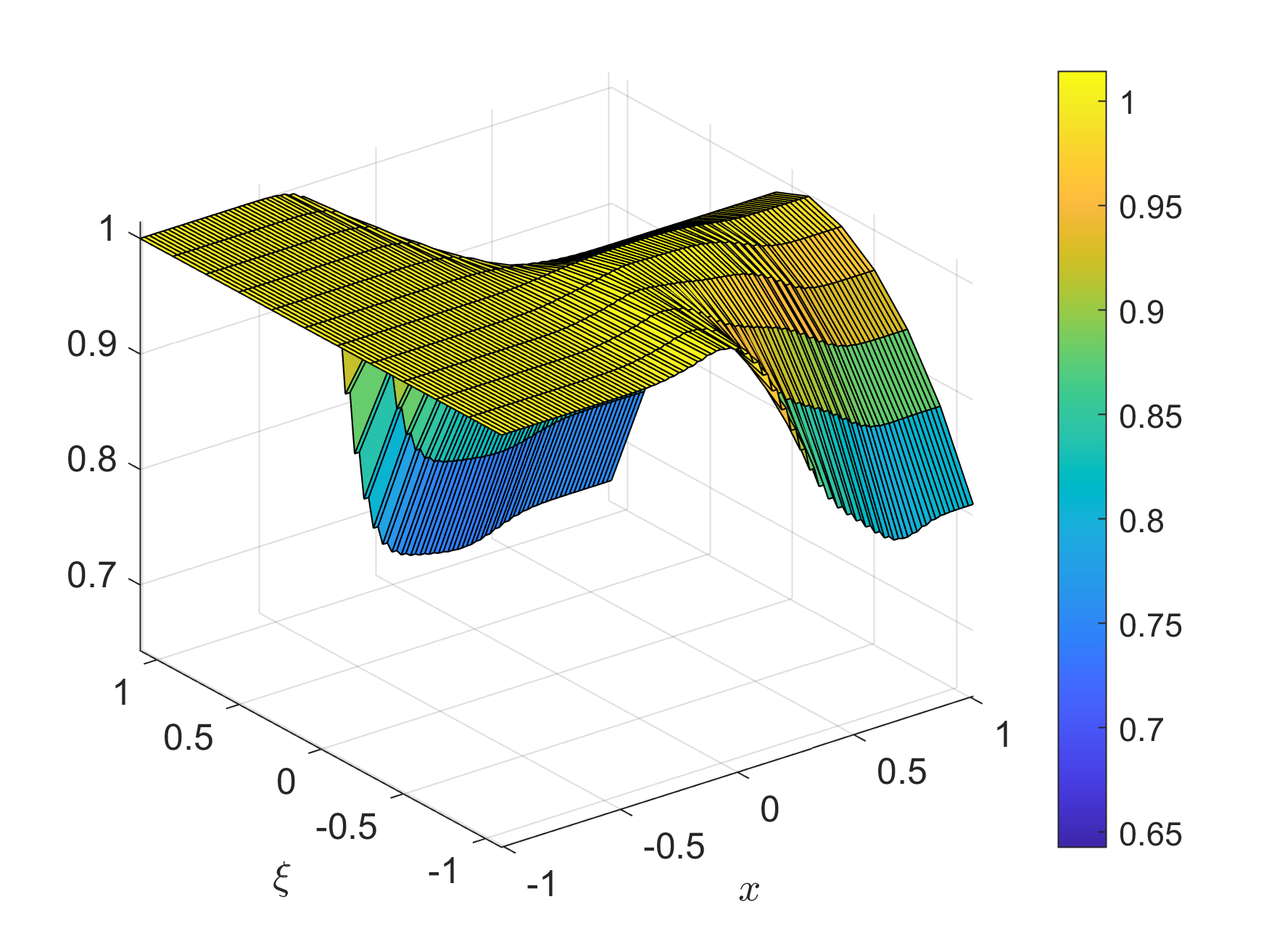}}
\caption{\sf Simulation results using the proposed method for the Euler system with the initial data given by \eqref{p-sys_id}. The solution $u(t,\cdot,\cdot)$ is shown at time $T=1/4.$ In the top part the initial data are shown. The second row presents the results of the stochastic collocation Lax-Friedrichs method and the bottom part shows the results obtained by the proposed method. To the left the density is shown and to the right the momentum is presented.\label{fig3}}
\end{figure}

The phase space for $u=(\rho,q)^T$ is chosen to be $[0.05, 2.5] \times [-1.0, 1.5]$ and discretized using $N_u=25$ points in {\em each} direction.  Terminal time is $T=1/4.$ The initial data for $(\rho,q)(0,\cdot,\cdot)$ are chosen to have 1-shocks and 1-rarefaction waves depending on the value of $\xi.$ To this end, we state the $1$-forward Lax curve starting at non-vacuum datum $U_L=(\rho_L,q_L).$ 
\begin{equation}\label{1-Lax}
\begin{aligned}
U(s; U_L) :=& \left( s, s \; \rho_L - \sqrt{ \frac{s}{\rho_L} ( s - \rho_L) ( p(s)-p(\rho_L) ) } \right),    \mbox{ if } s \geq \rho_L, \\
U(s; U_L) :=& \left( s, s \; \rho_L - s  \left(  \log(s) - \log(\rho_L)   \right) \right),  \mbox{ if } 0 < s < \rho_L. 
\end{aligned}
\end{equation}
Using formula \eqref{1-Lax} the following Riemann problem is considered. The left initial datum is  deterministic  $U_L=(1,1)^T$, and the right datum depends on the value of $\xi$. 
\begin{equation}\label{p-sys_id}
\begin{aligned}
u_0(x,\xi) :=& U_L, \; \mbox{ for } x<0 \mbox{ and all } \xi, \\
u_0(x,\xi) :=&  U( \frac12 \xi  + U_L; U_L), \mbox{ for } x>0 \mbox{ and all } \xi,
\end{aligned}
\end{equation}
The solution exhibits a shock wave for $\xi<0$ in the first family and a rarefaction wave of the first family  for all values of  $\xi\geq 0.$ In Figure \ref{fig3}, simulation results on the given spatial and random grids for both the proposed method as well as the reference stochastic collocation Lax-Friedrichs finite-volume method \eqref{coll} are shown. The time step $\Delta t$ is chosen according to the CFL condition \eqref{CFL2} and computed at each discrete point in time $t^n=n\Delta t.$ The $L^1$-error \eqref{L1err} in the density $\rho$ at time $T=1/4$ is $err=4.3707e$-$04$ and in the momentum $q$ it is $err=2.5379e$-$04$. The solution structure for different values of $\xi$ is also correctly identified with the proposed method.

\subsection{Conservation Laws with Discontinuous Flux}
It is instructive to examine whether the proposed scheme can capture non-unique numerical solutions in examples of conservation laws with discontinuous coefficients when different entropy functions are employed in the linear programming formulation. To this end, we consider the scalar conservation law with a discontinuous flux function:
$$
u_t+F(u)_x=0,
$$
where $F(u)=(1-H(x))g(u)+H(x)f(u)$, and $H(x)$ is the Heaviside function. We refer the reader to \cite{Mishra2017} for more details.

In this example taken from \cite{Mishra2017}, we take 
$$
g(u)=u(1-u), \quad  f(u)=1.1u(1-u),
$$ 
and consider the following initial data:
\begin{equation*}
 u_0(x) =
\begin{cases}
  0.65, & \mbox{if } x<0, \\
  0.35, & \mbox{otherwise},
\end{cases}
\end{equation*}
prescribed in the computational domain $[-4,4]$ subject to the free boundary conditions. The corresponding entropy functions $\eta$ in (\ref{linprog-aa}) are $\eta(u)=\frac{1}{2}u^2$ and $\eta(u)=|u-c|$, where $c$ is a constant satisfying $c\in [u_{\min}, u_{\max}]$.

\begin{figure}[ht!]
\centerline{\includegraphics[trim=0.2cm 0.2cm 0.8cm 0.1cm, clip, width=7.cm]{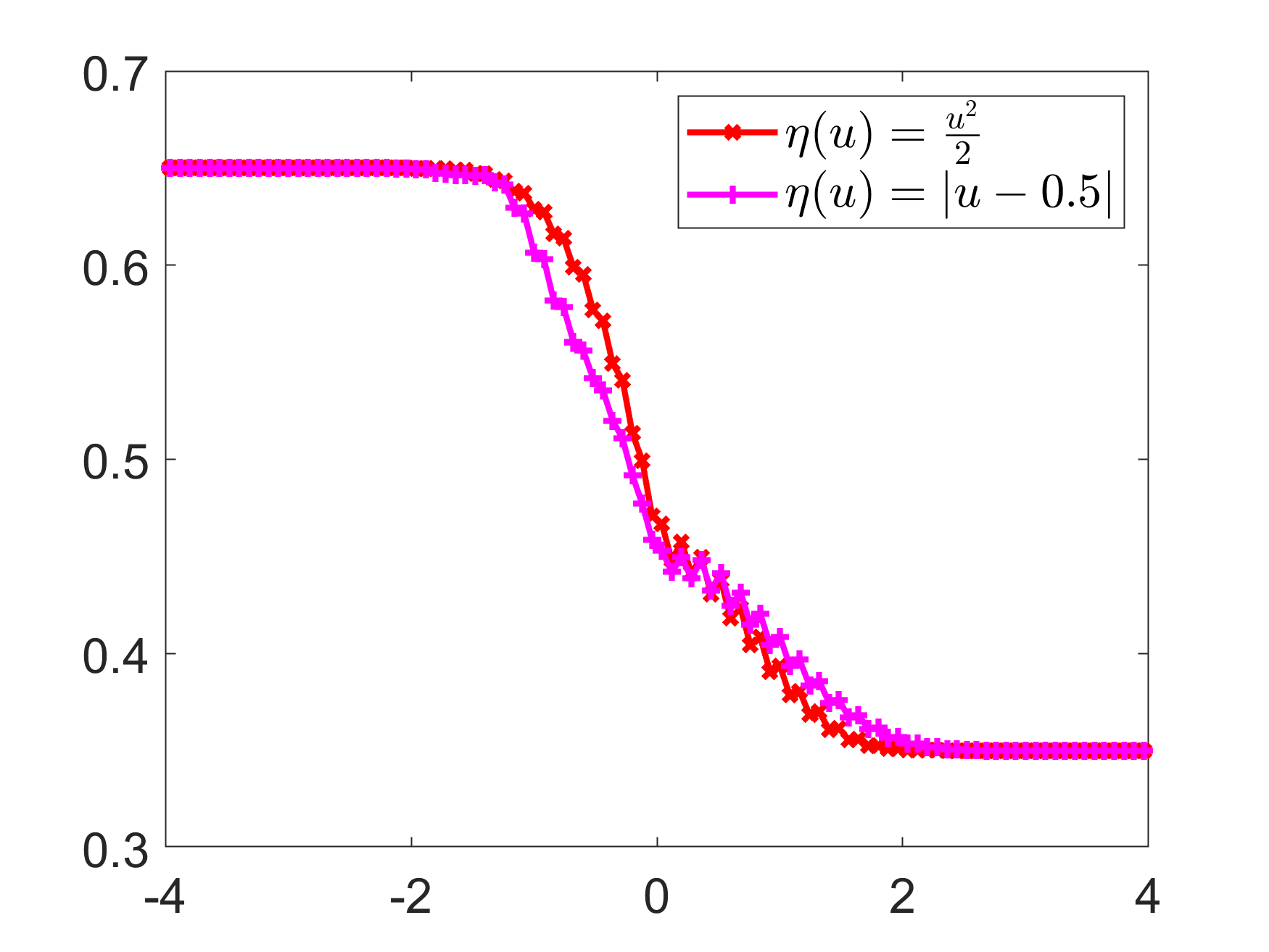}\hspace{0.5cm}
            \includegraphics[trim=0.2cm 0.7cm 0.8cm 0.1cm, clip, width=7.cm]{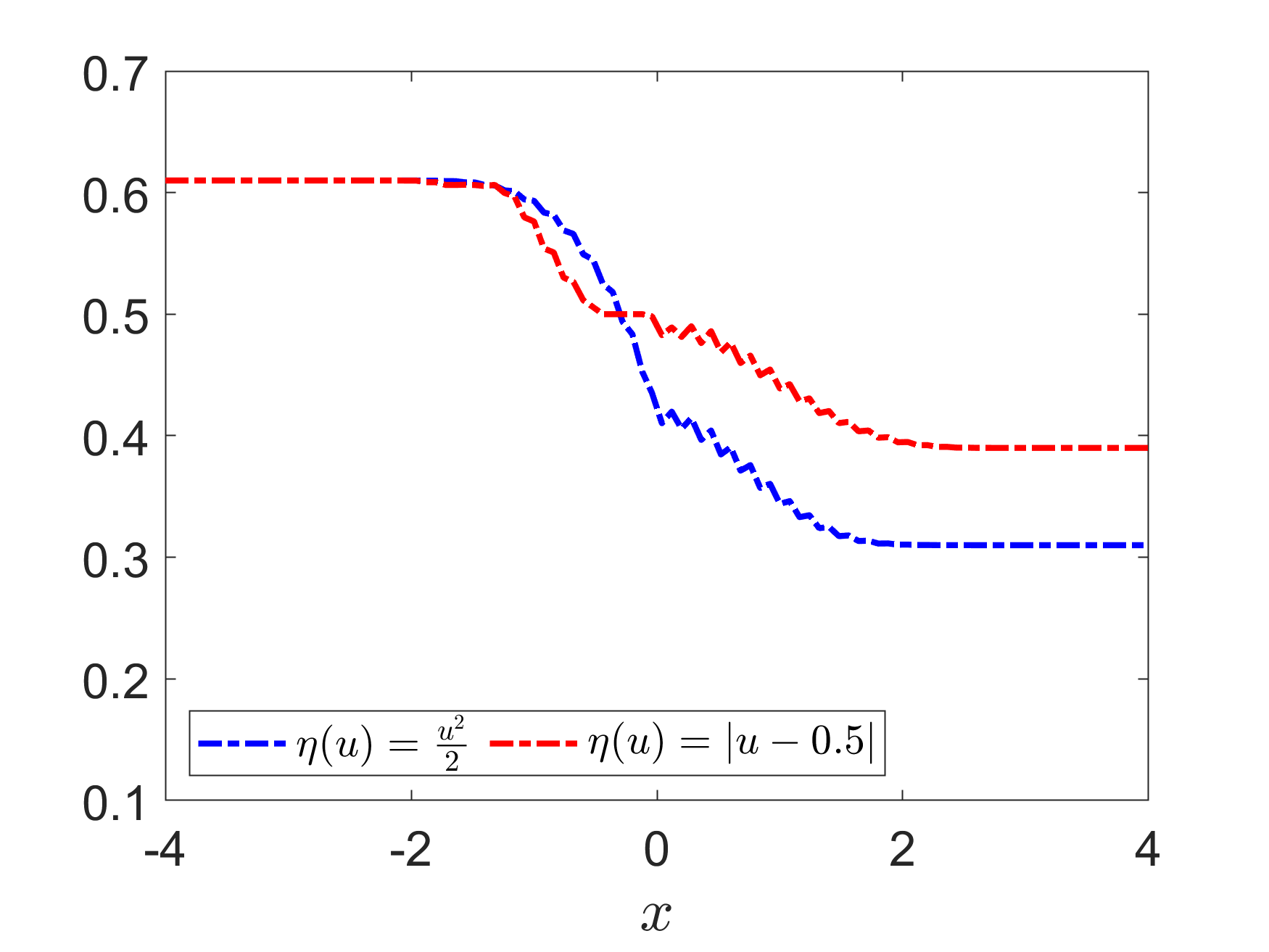}}
\caption{\sf Numerical results (left) and moments (right) at the final time $T=2$ with $\eta(u)=\frac{1}{2}u^2$ and $\eta(u)=|u-0.5|$. \label{fig3.11a}}
\end{figure}

\begin{figure}[ht!]
\centerline{\includegraphics[trim=0.2cm 0.2cm 0.8cm 0.1cm, clip, width=7.cm]{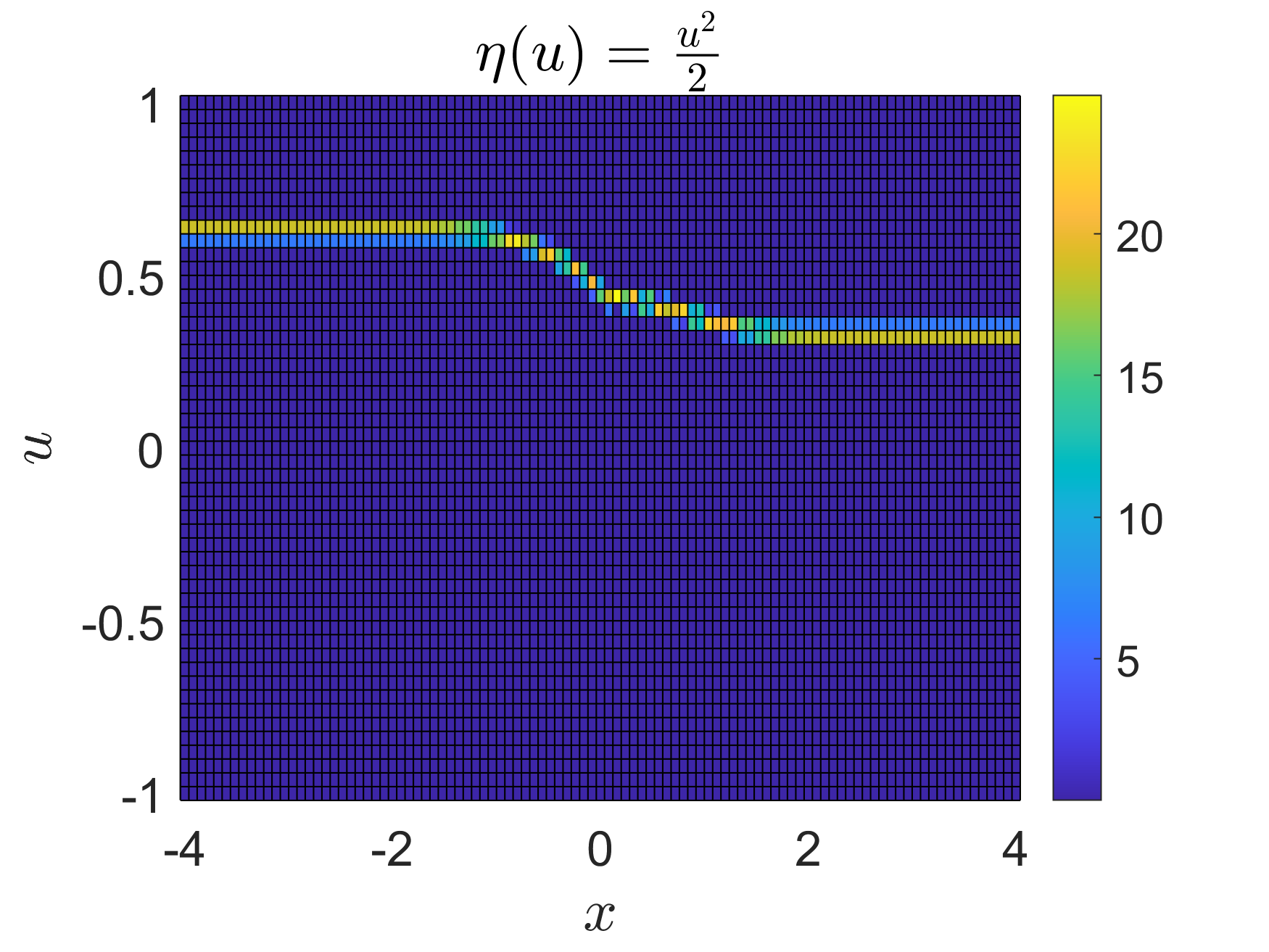}\hspace*{0.5cm}
            \includegraphics[trim=0.2cm 0.2cm 0.8cm 0.1cm, clip, width=7.cm]{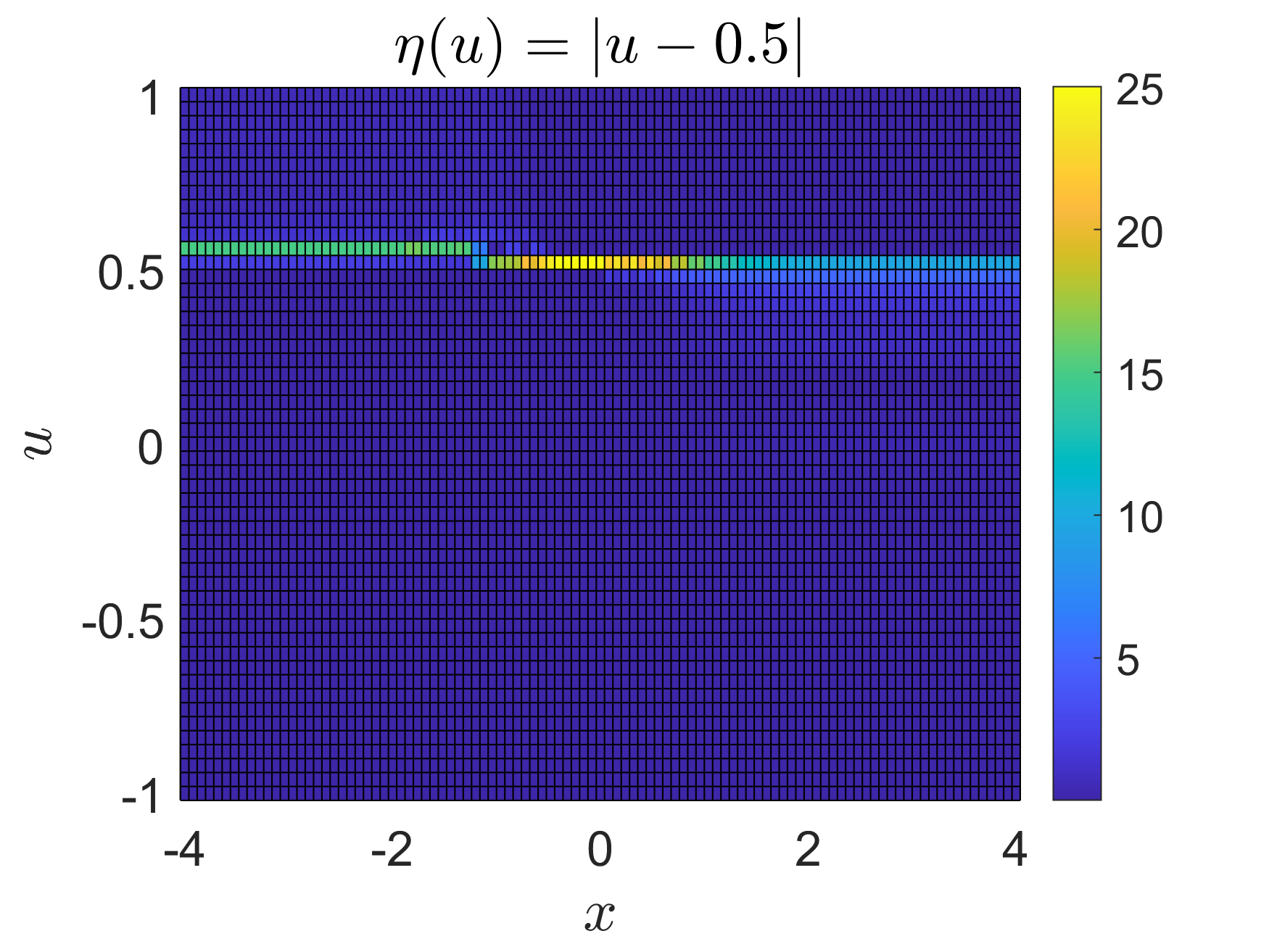}}
\caption{\sf Young measure at the final time $T=2$ with $\eta(u)=\frac{1}{2}u^2$ (left) and $\eta(u)=|u-0.5|$ (right). \label{fig3.11c}}
\end{figure}

We compute the numerical solutions until the final time $T=2$ on the uniform mesh with $N_x=100$. For the discretization in phase space, we set $u\in[-1,1]$ and use $N_{u}=50$. The obtained numerical results are presented in Figures \ref{fig3.11a}--\ref{fig3.11c}, where one can see that the numerical results are different when utilizing different entropy functions. The obtained results are in a good agreement with the numerical results reported in \cite{Mishra2017}, where varied schemes are used to produce different numerical solutions. We also present the corresponding moments and supports of the Young measures. One can notice that the moments and supports of the Young measures also vary when different entropy functions are utilized.

\begin{remark}
In the case of discontinuous fluxes, the  distinct solutions are obtained  by changing the objective function in the linear programming formulation. This example demonstrates that the 
proposed  formulation can serve as a selection criterion for the corresponding solutions.
\end{remark}

\section{Discussion and Outlook}
We have proposed a novel numerical method for random hyperbolic conservation laws based on the concept of measure-valued solutions. The method was presented for one-dimensional hyperbolic system and for one-dimensional random data. Generalization to multi-dimensional random space and multi-dimensional hyperbolic systems seems straightforward. In the context of the linear programs this would however lead to possibly very high dimensional problems. Generalization to high-order numerical methods is possible as well. The high-order reconstruction will enter into the moments that in turn are used in the constraints of the linear program. Any high-order numerical flux can then be computed using the solution to the linear program.

Furthermore, in this work only numerical results for the basis functions $\phi_i(\xi) =\delta(\xi-\xi_i)$ (see \eqref{base-coll}) for its discrete representation, have been shown. However, our approach is more general and can use any orthonormal basis in random space. Hence, the suggested method yields a unifying framework for the treatment of random hyperbolic conservation laws. In particular, the proposed method is less intrusive than the stochastic Galerkin method using generalized polynomial chaos, since only the flux function needs to be replaced by the linear programming problem.  The conservative variables and their high-order reconstruction do not  require any modification. Of course, it is still more intrusive than the Monte Carlo or stochastic collocation methods. The numerical results show that in the considered case of basis functions the numerical solutions obtained by the new method and the standard stochastic collocation method are comparable. The proposed method is computationally less efficient due to the subsequent solutions of possibly large linear programs. 
 
We point out that our technique uses the underlying measure-valued representation of a solution and approximates the corresponding Young measure directly. Here, we proposed a direct computation of the Young measure using a particular discretization of the measure. There are also approaches available in the literature as outlined in the introduction which compute the measure using a sequence of moments as well as simulation on different spatial grids. Future work will assess the difference between these approaches and the proposed method.
 
{\small \subsection*{Acknowledgments}
The authors thank the Deutsche Forschungsgemeinschaft (DFG, German Research Foundation) for the financial support through 320021702/GRK2326,  333849990/GRK2379, SPP 2410 Hyperbolic Balance Laws in Fluid Mechanics: Complexity, Scales, Randomness (CoScaRa) within the Projects HE5386/26-1 (Numerische Verfahren f\"ur gekoppelte Mehrskalenprobleme, 525842915),  HE5386/27-1, LU1470/10-1  (Zuf\"allige kompressible Euler Gleichungen: Numerik und ihre Analysis, 525853336)  and LU1470/9-1 (An Active Flux method for the Euler equations, 525800857). M.H. has also received funding from the European Union's Horizon Europe research and innovation programme under the Marie Sklodowska-Curie Doctoral Network Datahyking (Grant No.101072546) and HE5386/30-1  Effiziente Ermittlung von Temperaturfeldern mithilfe von Physics-Informed Neural Networks zur Modellierung des thermo-elastischen Verhaltens von Werkzeugmaschinen (537928890). 
The work  of M.L.-M. was also partially supported by the Gutenberg Research College and by the Deutsche Forschungsgemeinschaft - project number 233630050 - TRR 146.  M.L.-M. is grateful to  the  Mainz Institute of Multiscale Modelling  for supporting her research.
}

\bibliographystyle{siam}
\bibliography{references}
\end{document}